\documentclass[11pt, twoside, leqno]{article}
\usepackage{amssymb,latexsym}
\usepackage{enumerate}
\usepackage{color}
\usepackage{amsmath, amsthm, amsfonts, amssymb, mathrsfs}
\newtheorem{theorem}{Theorem}[section]

\newtheorem{lemma}[theorem]{Lemma}
\newtheorem{proposition}[theorem]{Proposition}

\newtheorem{remark}[theorem]{Remark}
\newtheorem*{thma}{\hskip\parindent {Theorem A}}

\def\dint{\displaystyle\int}
\def\dsum{\displaystyle\sum}
\def\dsup{\displaystyle\sup}

\numberwithin{equation}{section}

\numberwithin{equation}{section}

\numberwithin{equation}{section}

\begin{document}
%\\\\\\\\\\\\\\title\\\\\\\\\\\\\\
\title{{\LARGE\bf    Weighted jump and variational inequalities for rough operators
 \footnotetext{\small
{{\it MR(2000) Subject Classification}:\ } 42B20, 42B25}
\thanks {The research was supported by
NSF of China (Grant: 11471033, 11371057, 11571160, 11401175, 11601396, 11431011, 11501169), NCET of China (Grant: NCET-11-0574),  the Fundamental Research Funds for the Central Universities (FRF-BR-16-011A, 2014KJJCA10), SRFDP of China   (Grant: 20130003110003), Thousand Youth Talents Plan  of China (Grant: 429900018-101150(2016)),  Funds for Talents of China (Grant: 413100002).}}}

\date{Revised Version }
\maketitle

\begin{center}
{\bf Yanping Chen \footnote{\small {Corresponding author.\ }}}\\
Department of Applied Mathematics,  School of Mathematics and Physics,\\
 University of Science and Technology Beijing,\\
Beijing 100083, The People's Republic of China \\
 E-mail: {\it yanpingch@ustb.edu.cn}
\vskip 0.1cm

{\bf Yong Ding}\\
School of Mathematical Sciences, Beijing Normal University,\\
Laboratory of Mathematics and Complex Systems (BNU),
Ministry of Education,\\
Beijing 100875,  The People's Republic of China \\
E-mail: {\it dingy@bnu.edu.cn}\vskip 0.1cm

{\bf Guixiang Hong}\\
School of Mathematics and Statistics, Wuhan University,\\
Wuhan 430072,  The People's Republic of China \\
E-mail: {\it guixiang.hong@whu.edu.cn}\vskip 0.1cm
 and \vskip 0.1cm

{\bf Honghai Liu}\\
School of Mathematics and Information Science,
Henan Polytechnic University,\\
Jiaozuo, Henan, 454003,  China \\
E-mail: {\it hhliu@hpu.edu.cn}\end{center}

\begin{center}
\begin{minipage}{135mm}{\small
{\bf ABSTRACT\ }In this paper, we systematically study  weighted jump and variational inequalities for rough operators. More precisely, we show some weighted jump and variational inequalities for the families $\mathcal T:=\{T_\varepsilon\}_{\varepsilon>0}$ of truncated singular integrals and $\mathcal M_{\Omega}:=\{M_{\Omega,t}\}_{t>0}$ of averaging operators with rough kernels, which are defined  respectively by
$$
T_\varepsilon
f(x)=\int_{|y|>\varepsilon}\frac{\Omega(y')}{|y|^n}f(x-y)dy
$$
 and
$$M_{\Omega,t} f(x)=\frac1{t^n}\int_{|y|<t}\Omega(y')f(x-y)dy,
$$
where the kernel $\Omega$ belongs to $L^q(\mathbf S^{n-1})$ for $q>1$.  }
\end{minipage}
\end{center}
\vspace{0.3cm}

%%%%%%%%%%%%%%%%%%%%%%%%%%%%  ?????? ????%%%%%%%%%%%%%%%%%%%%%
%%%%%%%%%%%%%%%%%%%%%%%%%%%%  ?????? ????%%%%%%%%%%%%%%%%%%%%%

 \pagestyle{plain} %???????
 \renewcommand{\thepage}{\arabic{page}}% ?????????????????
 \setcounter{page}{1} %?????1???
 \setcounter{equation}{0}%?????1???????

 \section{Introduction}\label{sect1}

The jump and variational inequalities have been the subject of many
 recent articles in probability, ergodic theory and harmonic
 analysis. To present related results in a precise way, let us fix some notations.
 Given a family of complex numbers $\mathfrak{a}=\{a_t: t\in\mathbb{R}\}$ and $\rho\ge1$, the $\rho$-variation norm of the family $\mathfrak{a}$ is defined by
\begin{equation}\label{q-ver number family}
\|\mathfrak{a}\|_{V_\rho}=\sup\big(|a_{t_0}|+\sum_{k\geq1}
|a_{t_k}-a_{t_{k-1}}|^\rho\big)^{\frac{1}{\rho}},
\end{equation}
where the supremum runs over all finite increasing sequences $\{t_k:k\geq0\}$. It is trivial that
\begin{equation}\label{number contr ineq}
\|\mathfrak{a}\|_{L^\infty(\mathbb{R})}:=\sup_{t\in\mathbb R}|a_t|
\le\|\mathfrak{a}\|_{V_\rho}\quad\text{for}\ \ \rho\ge1.
\end{equation}

Via the definition \eqref{q-ver number family}  of the $\rho$-variation norm of a family of numbers, one may define the strong $\rho$-variation function $V_\rho(\mathcal F)$ of a family $\mathcal F$ of functions.
Given a family of Lebesgue measurable functions $\mathcal F=\{F_t:t\in\mathbb{R}\}$ defined on $\mathbb{R}^n$, for any fixed $x$ in $\mathbb{R}^n$, the value of the strong $\rho$-variation function $V_\rho(\mathcal F)$ of the family $\mathcal F$ at $x$ is defined by
\begin{equation}\label{defini of Vq(F)}
V_\rho(\mathcal F)(x)=\|\{F_t(x)\}_{t\in\mathbb{R}}\|_{V_\rho},\quad \rho\ge1.
\end{equation}

 Suppose $\mathcal{A}=\{{A}_t\}_{t>0}$ is a family of operators  on $L^p(\Bbb R^n)\, (1\le p\le\infty)$. The strong $\rho$-variation operator is simply defined as
$$V_\rho(\mathcal Af)(x)=\|\{A_t(f)(x)\}_{t>0}\|_{V_\rho},\quad\forall f\in L^p(\mathbb{R}^n).$$
It is easy to observe from the definition of $\rho$-variation norm that for any $x$ if $V_\rho(\mathcal Af)(x)<\infty$, then $\{A_t(f)(x)\}_{t>0}$ converges when $t\rightarrow0$ or $t\rightarrow\infty$. In particular, if $V_\rho(\mathcal Af)$ belongs to some function spaces such as $L^p(\mathbb{R}^n)$ or $L^{p,\infty}(\mathbb{R}^n)$, then the sequence converges almost everywhere without any additional condition. This is why the mapping property of the strong $\rho$-variation operator is so interesting in ergodic theory and harmonic analysis. Also, by \eqref{number contr ineq}, for any $f\in L^p(\Bbb R^n)$ and $x\in\Bbb R^n$, we have
 \begin{equation}\label{control of maxi opera}
A^\ast(f)(x)\le V_{\rho}(\mathcal{A}f)(x)\quad\text{for}\ \ \rho\ge1,
\end{equation}
where $A^\ast$ is the maximal operator defined by
$$A^\ast(f)(x):=\sup_{t>0}|{A}_t(f)(x)|.$$

Let $\mathcal F=\{F_t(x):t\in\mathbb{R}_+\}$ be a family of Lebesgue measurable functions defined on $\mathbb{R}^n$. For $\lambda>0$, we introduce the $\lambda$-jump function $N_\lambda(\mathcal F)$ of $\mathcal F$, its value at $x$ is the supremum over all $N$ such that there exist $s_1<t_1\leq s_2<t_2\leq\dotsc\leq s_N<t_N$ with
$$|F_{t_k}(x)-F_{s_k}(x)|>\lambda$$
for all $k=1,\dotsc, N$.

Using the fact that $\ell^{2,\infty}(\mathbb{N})$ (the weak $L^2$ space on $\Bbb N$) embeds into $\ell^\rho(\mathbb{N})$ for all $\rho>2$, it is easy to check that the following pointwise domination holds
$$\sup_{\lambda>0}\lambda \sqrt{N_\lambda(\mathcal F)}\geq V_\rho(\mathcal F)$$
for all  families of functions $\mathcal F$.

 The first variational inequality was proved by L\'epingle \cite{Lep76} for martingales (see \cite{PiXu88} for a simple proof). Bourgain \cite{Bou89} is the first one using L\'epingle's result to obtain similar variational estimates for the ergodic averages,  and then directly deduce pointwise convergence results without previous knowledge that the pointwise convergence holds for a dense subclass of functions, which are not available in some ergodic models. In particular, Bourgain's work \cite{Bou89} has inaugurated a new research direction in ergodic theory and harmonic
analysis. In their papers \cite{JKRW98, JRW03, JRW00, CJRW2000, CJRW2002}, Jones and his collaborators systematically studied jump and variational inequalities for ergodic averages and truncated singular integrals (mainly of homogeneous type). Since then many other publications came to enrich the literature on this subject (cf. e.g. {\cite{GiTo04, LeXu2, DMT12,JSW08, MaTo, OSTTW12,  MaS, MaBP}).

For $\varepsilon>0,$ suppose $T_\varepsilon$ is the truncated singular integral operator defined by
 \begin{equation}\label{tr of S}
T_\varepsilon
f(x)=\int_{|y|>\varepsilon}\frac{\Omega(y')}{|y|^n}f(x-y)dy,
 \end{equation}
where $\Omega\in L^1({\mathbf S}^{n-1})$
 satisfies the cancelation condition
\begin{equation}\label{can of O}
\int_{\mathbf S^{n-1}}\Omega(y')d\sigma(y')=0.
 \end{equation}
 Denote the family of operators $\{T_\varepsilon\}_{\varepsilon>0}$ by $\mathcal T$. For $1<p<\infty$ and $f\in L^p(\mathbb R^n)$,
 the Calder\'on-Zygmund singular integral operator $T$ with homogeneous kernel is defined by
\begin{equation}\label{SIO}
Tf(x)=\lim_{\varepsilon\rightarrow0^+}T_\varepsilon f(x), \ a.e. \ x\in\mathbb R^n.
\end{equation}
The fact that the limit exists almost everywhere is deduced from the density of $C^\infty_0(\mathbb R^n)$ in $L^p(\mathbb R^n)$ and the boundedness of the maximal singular integral $T^\ast$ which is given by $T^\ast f(x)=\sup_{\varepsilon>0}|T_\varepsilon f(x)|$. In fact, variational and jump inequalities for $\mathcal T$ can also be used to study the existence of the above limit and give extra information on the convergence.

The famous Hilbert transform $H$, which is defined by
$
H(f)(x)={\rm p.v.}\frac 1\pi\int_{\mathbb{R}}\frac{f(y)}{x-y}dy,
$
is an example of  homogeneous singular integral operators $T_\Omega$ when the dimension $n=1$. In 2000, Campbell {\it et al} \cite{CJRW2000} first considered the $L^p(\mathbb R)\,(1<p<\infty)$ boundedness of the strong $\rho$-variation operator of the family of the truncated Hilbert transforms denoted by $\mathcal H:=\{H_\varepsilon\}_{\varepsilon>0}$.
In 2002, Campbell {\it et al} \cite{CJRW2002} gave the $L^p(\mathbb R^n)$ boundedness of the strong $\rho$-variation operator of  $\mathcal T$, the family of homogenous singular integrals with $\Omega\in L\log^+\!\!L(\mathbf S^{n-1})$ and $n\geq2$ for $\rho>2$.
In 2008, using the Fourier transform and the square function estimates given in \cite{DR86},
 Jones, Seeger and Wright \cite{JSW08} developed a general method, which allows one to obtain some jump inequalities for  families of the truncated singular integral operators $\mathcal T$ and of other integral operators arising from harmonic analysis.

\begin{thma} {\rm (\cite{JSW08})}
Suppose $\Omega$ satisfies \eqref{can of O} and $\Omega\in L^r(\mathbf S^{n-1})$ for $r>1$. Then the $\lambda$-jump inequality  $$\sup_{\lambda>0}\|\lambda \sqrt{N_\lambda(\mathcal T f)}\|_{L^p(\mathbb R^n)}\le
C_{p,n}\|f\|_{L^p(\mathbb R^n)}\ (1<p<\infty)$$ holds,
where and in the sequel, the constant $C_{p,n}>0$ depends only on $p$ and $n$. Whence, for $\rho>2$ and $1<p<\infty$, there
exists a constant $C(p,\rho)$ such that
\begin{equation*}
 \|V_\rho(\mathcal Tf)\|_{L^p(\mathbb R^n)}\le
C(p,\rho)\|f\|_{L^p(\mathbb R^n)}.
\end{equation*}
\end{thma}

The purpose of this paper is to give the  weighted jump and variational inequalities for singular integrals and averaging operators with rough kernels.
In order to state our main results, let us first recall some definitions. We first recall the definition and some properties of $A_p$ weight on $\mathbb R^n$. Let $w$ be a non-negative locally integrable function defined on $\mathbb R^n$. We say $w\in A_1$ if there is a constant $C>0$ such that $M(w)(x)\le Cw(x)$, where $M$ is the classical Hardy-Littlewood maximal operator defined by
$$Mf(x)=\sup_{r>0}\frac 1{r^n}\int_{|y|\le r}|f(x-y)|dy.$$
 Equivalently, $w\in A_1$ if and only if there is a constant $C>0$ such that for any cube $Q$
\begin{equation}\label{eq def A1}
\frac1{|Q|}\int_Qw(x)dx\le C\inf_{x\in Q}w(x).
\end{equation}
For $1<p<\infty$, we say that $w\in A_p$ if there exists a constant $C>0$ such that
\begin{equation}\label{eq def Ap}
\sup_Q\bigg(\frac1{|Q|}\int_Qw(x)dx\bigg)\bigg(\frac1{|Q|}\int_Qw(x)^{1-p'}dx\bigg)^{p-1}\le C.
\end{equation}
The smallest constant appearing in \eqref{eq def A1} or \eqref{eq def Ap} is denoted by $[w]_{A_p}$.
$A_\infty=\bigcup_{p\ge1} A_p$. It is well known that if $w\in A_\infty$, then there exist $\delta\in(0,1]$ and $C>0$ such that for any cube $Q$ and measurable subset $E\subset Q$
\begin{equation}\label{pap}
\frac{w(E)}{w(Q)}\le C\bigg(\frac{|E|}{|Q|}\bigg)^\delta.
\end{equation}

The boundedness of many important operators in harmonic analysis on $L^p(w)
$ for $w\in A_p$, $1 <p<\infty$, has been known for a long time. We
refer the reader to \cite{Gar, GR, MW, Mu1, Mu2, D93, Wa, GMS} for more details on this topic. Recently there has been renewed interest in  weighted  jump and variational inequalities.
The  weighted $\rho$-variational inequality $2<\rho<\infty$ for singular integrals with Lipschitz kernels has been shown recently in \cite{MTX1, MTX2} (see also \cite{KZ14} for averaging operators). It arises naturally as an open problem that whether the $\rho$-variations for singular integrals with rough kernels are bounded on weighted $L^p$ spaces, even though the boundedness on unweighted $L^p$ spaces has been proved in \cite{CJRW2000, CJRW2002, JSW08, DHL}. In the present paper, we give a positive solution to the problem  and  prove a weighted jump inequality which implies all  $\rho$-variational inequalities for the singular integrals with rough kernels. We also show  corresponding weighted jump inequalities for the  averaging operators with rough kernels.

\begin{theorem}\label{thm:H}  Let $\mathcal T$ be given as
in \eqref{tr of S} with $\Omega\in L^q(\mathbf S^{n-1})$, $q>1$ satisfying \eqref{can of O}.
Then the following $\lambda$-jump inequality holds
\begin{equation}\label{N12w}
\|\sup_{\lambda>0}\lambda \sqrt{N_\lambda(\mathcal T f)}\|_{L^p(w)}\le
 C_{p,w}\|f\|_{L^p(w)},
\end{equation}
if $w$ and $p$ satisfy one of the following conditions:

{\rm(i)} If $q'\le p<\infty,\, p\neq 1$  and $ w\in A_{p/q'},$

{\rm (ii)} If $1< p\le q,\, p\neq \infty$  and $ w^{-\frac{1}{(p-1)}}\in A_{p'/q'}$.

Whence, for $\rho>2$, there
exists a constant $C(p,\rho)$ such that
\begin{equation}\label{Vqr}
 \|V_\rho(\mathcal Tf)\|_{L^p(w)}\le
C(p,\rho)\|f\|_{L^p(w)},
\end{equation} if $w$ and $p$ satisfy one of the conditions (i) or (ii).
\end{theorem}

\begin{remark} Theorem \ref{thm:H} covers Theorem A with $w\equiv1$.  Also, \eqref{Vqr} is an improvement of the weighted $L^p$ boundedness for $T^\ast$(see \cite{D93} and \cite{Wa}), because of the pointwise estimate $T^\ast f(x)\le V_\rho(\mathcal Tf)(x)$.  Restricted to the singular integrals of homogeneous type, this result significantly improves Corollary 1.4 in  \cite{MTX2} where $\Omega$ is assumed to be in the H\"older class of order $\alpha$.\end{remark}

As in \cite{JSW08}, the first step to prove Theorem \ref{thm:H} is that the desired estimate (\ref{N12w}) is reduced to the estimate over short $2$-variation and the estimate over dyadic $\lambda$-jump function through the following pointwise inequality (see for instance Lemma 1.3 in \cite{JSW08})
\begin{equation}\label{lem:convert lemma}
\lambda\sqrt{N_\lambda(\mathcal Tf)(x)}\le
C\big[S_2(\mathcal Tf)(x)+\lambda\sqrt{N_{\lambda/3}(\{
T_{2^k}f\}_{k\in \Bbb Z})(x)}\big],
 \end{equation}
where
 $$
S_2(\mathcal Tf)(x)=\bigg(\sum_{j\in\mathbb Z}[V_{2,j}(\mathcal
Tf)(x)]^2\bigg)^{1/2},
$$
with
$$
V_{2,j}(\mathcal Tf)(x)=\bigg(\sup_{\substack
{t_1<\cdots<t_N\\
[t_l,t_{l+1}]\subset[2^j,2^{j+1}]}}\sum_{l=1}^{N-1}|T_{t_{l+1}}f(x)-T_{t_l}f(x)|^2\bigg)^{1/2}.
 $$
That is, we are reduced to prove
\begin{align}\label{dyadic jump singular}
\|\sup_{\lambda>0}\lambda \sqrt{N_\lambda(\{T_{2^k}f\})}\|_{L^p(w)}\le
C_{p,w}\|f\|_{L^p(w)}
\end{align}
and
\begin{align}\label{short variation singular}
\|S_2(\mathcal Tf)\|_{L^p(w)}\le
 C_{p,w}\|f\|_{L^p(w)}
\end{align}
for $w$ and $p$ satisfying the conditions (i) or (ii) in Theorem \ref{thm:H}.

 In order to show estimate (\ref{dyadic jump singular}), we prove some vector-valued weighted estimates such as (\ref{Tskfk}) in Section 2 and
use the generalized Rubio de Francia's extrapolation theorem---Lemma \ref{extrapolation} in Section 2--- as well as Stein and Weiss's interpolation theorem with change of measure. On the other hand, to establish estimate (\ref{short variation singular}), we discover a new phenomenon, that is, the short $2$-variation can be dominated by the vector-valued maximal function with rough kernel---inequality (\ref{s2km}) in Section 3, which actually provides an alternate proof instead of the rotation argument used in \cite{CJRW2002} and \cite{JSW08}. Then we are reduced to prove vector-valued weighted estimates such as (\ref{m2pw}) and (\ref{Tjt2pw}) in Section 3. The main idea behind the proof is Rubio de Francia's extrapolation theorem, see for instance Remark \ref{one key remark} below.

\bigskip

The proof of Theorem \ref{thm:H} can be adapted to the situation of averaging operators with rough kernels $\mathcal M_\Omega=\{M_{\Omega,t}\}_{t>0}$, where $M_{\Omega,t}$ is defined as
\begin{equation}\label{def M O}
M_{\Omega,t} f(x)=\frac1{t^n}\int_{|y|<t}\Omega(y')f(x-y)dy,
\end{equation}
where $\Omega\in L^1({\mathbf S}^{n-1})$.

\begin{theorem}\label{r o MO}
 Suppose the family $\mathcal M_\Omega=\{M_{\Omega,t}\}_{t>0}$ is defined in \eqref{def M O}.
Let $\Omega\in L^q(\mathbf S^{n-1}) $ for $q>1.$  Then
\begin{equation}
\|\sup_{\lambda>0}\lambda \sqrt{N_\lambda(\mathcal M_\Omega f)}\|_{L^p(w)}\le
 C_{p,w}\|f\|_{L^p(w)}
\end{equation}
if $w$ and $p$ satisfy (i) or (ii) in Theorem \ref{thm:H}.
The similar inequality holds for the strong $\rho$-variation operator $V_\rho(\mathcal M_\Omega f)$ with $\rho>2$.
\end{theorem}

\begin{remark} When $\Omega\equiv1$, the weighted jump and variational inequalities have been proved in \cite{MTX1}, \cite{MTX2} and \cite{KZ14}. We will explain briefly the proof of Theorem \ref{r o MO} in Section 4.\end{remark}

\section{Proof of Theorem \ref{thm:H} (I)}

As we have stated in the previous section, to prove Theorem \ref{thm:H} it suffices to show \eqref{dyadic jump singular} and \eqref{short variation singular}. In this section, we give the proof of \eqref{dyadic jump singular}. Let us begin with one definition. For
$j\in\mathbb Z$,  let $\nu_j(x)=\frac{\Omega(y)}{|y|^n}\chi_{\{2^j\le |x|<2^{j+1}\}}(x)$, then
$$
\nu_j\ast f(x)=\int_{2^j\le |y|<2^{j+1}}\frac{\Omega(y)}{|y|^n}f(x-y)dy.
$$
Obviously, for $k\in\mathbb Z$,
$$
T_{2^k}f(x)=\int_{|x-y|\ge
2^k}\frac{\Omega(x-y)}{|x-y|^n}f(y)\,dy=\sum_{j\ge k}\nu_j\ast f(x).
$$
Let $\phi\in \mathscr{S}(\mathbb R^n)$ be a radial function such that $\hat{\phi}(\xi)=1$ for
$|\xi|\le2$ and $\hat{\phi}(\xi)=0$ for $|\xi|>4$. We have the following decomposition
\begin{align*}
T_{2^k}f&=\phi_k\ast Tf+\sum_{s\ge0}(\delta_0-\phi_k)\ast\nu_{k+s}\ast f-\phi_k\ast\sum_{s<0}\nu_{k+s}\ast
f\\
&:=T^1_{k}f+T^2_{k}f-T_{k}^3f,
\end{align*}
where $\phi_k$ satisfies $\widehat{\phi_k}(\xi)=\hat{\phi}(2^k\xi)$, $\delta_0$ is the Dirac measure at 0 and $s\in\Bbb N\cup \{0\}$. Let $\mathscr T^if$ denote the family $\{T^i_{k}f\}_{k\in\mathbb Z}$  for
$i=1,2,3$. Obviously, to show \eqref{dyadic jump singular} it suffices to prove the following inequalities:
\begin{equation}\label{lnti}
\|\sup_{\lambda>0}\lambda [N_{\lambda}(\mathscr T^if)]^{1/2}\|_{L^p(w)}\le C_{p,w}\|f\|_{L^p(w)},\ \ i=1,2,3,
\end{equation}
for $w$ and $p$ satisfying the conditions (i) or (ii) in Theorem \ref{thm:H}.

\bigskip

\noindent {\bf Estimate of \eqref{lnti} for $i=1$.}
This estimate will follow easily from the weighted $L^p$-boundedness of $T$ (see \cite{D93} or \cite{Wa} ) and the following Proposition \ref{pro:DNP} which is a simplified and weighted version of Theorem 1.1 in \cite{JSW08}. The proof of Proposition \ref{pro:DNP} will be postponed to the end of the section.
\begin{proposition}\label{pro:DNP}
Let $\mathscr U$ be a family of operators given by $\mathscr Uf=\{\phi_k\ast f\}_k$. Then for
$1<p<\infty$ and $w\in A_p$, we have
\begin{equation}\nonumber
\|\sup_{\lambda>0}\lambda\sqrt{N_\lambda(\mathscr{U}f)}\|_{L^p(w)}\leq C_{p,w}\|f\|_{L^p(w)}.
\end{equation}
\end{proposition}
Indeed,
\begin{align*}
\|\sup_{\lambda>0}\lambda[N_\lambda(\mathscr T^1f)]^{1/2}\|_{L^p(w)}&\le\|\sup_{\lambda>0}\lambda[N_\lambda(\{\phi_k\ast  Tf\})]^{1/2}\|_{L^p(w)}\\&\le  C_{p,w}\|  Tf\|_{L^p(w)}\le  C_{p,w}\|\Omega\|_{L^q(\mathbf S^{n-1})}\|f\|_{L^p(w)},
\end{align*}
 whenever $w$ and $p$ satisfy the conditions (i) or (ii) in Theorem \ref{thm:H}.

\medskip

\noindent {\bf Estimate of \eqref{lnti} for $i=2$.}
By the Minkowski inequality, we get
\begin{align}
\nonumber\sup_{\lambda>0}\lambda[N_\lambda(\mathscr T^{2}f)(x)]^{1/2}&\le\dsum_{s\ge0}\Big(\dsum_{k\in \Bbb Z}\Big|(\delta_0-\phi_k)\ast\nu_{k+s}\ast
f(x)\Big|^2\Big)^{1/2}\\
\label{lnt2}&:=\dsum_{s\ge0}G_sf(x).
\end{align}

We are reduced to establish  the estimate of $\|G_sf\|_{L^p(w)}$ as sharp as possible so that we are able to sum up over $s\in\Bbb N\cup \{0\}$.
Let $\psi\in
C_0^\infty({\Bbb R}^n)$ be a radial function such
  that
 $0\le \psi\le 1,$ $\hbox{supp}\psi \subset\{1/2\le |\xi|\le
2\}$ and $\sum_{l\in \Bbb Z}\psi^2(2^{-l}\xi)=1$ for $|\xi|\neq 0.$
Define the multiplier $\Delta_l$ by
$ \widehat{\Delta_l f}(\xi)=\psi(2^{-l}\xi)\widehat{f}(\xi)$. Therefore,
\begin{align}
\nonumber G_sf(x)&=\bigg(\sum_{k\in \Bbb Z}\big|[(\delta_0-\phi_k)\ast\nu_{s+k}]\ast
\sum_{l\in \Bbb Z}\Delta_{l-k}^2f(x)\big|^2\bigg)^{1/2}\\
\label{Gsl}&\le \sum_{l\in \Bbb Z}\bigg(\sum_{k\in \Bbb Z}|\Delta_{l-k}[(\delta_0-\phi_k)\ast\nu_{s+k}]\ast
\Delta_{l-k}f(x)|^2\bigg)^{1/2}:=\sum_{l\in \Bbb Z}G_s^{l}f(x).
\end{align}
 We first prove a rapid decay estimate of $\|G_{s}^{l}f\|_{L^2(\mathbb R^n)}$ for $l\in \Bbb Z$ and $s\in\Bbb N\cup \{0\}$. Since
 $$\hbox{supp}(1-\widehat{\phi_k})\widehat{\nu_{k+s}} \subset \{\xi :|2^k\xi|>2\}$$ and $\Omega(x')$ satisfies \eqref{can of O},
by a well-known Fourier transform estimate of Duoandikoetxea and Rubio de Francia
 (see \cite[p.551-552]{DR86}), it is easy to show that for any fixed $q>1$ and  some $\gamma\in(0,1)$,
$$|1-\widehat{\phi_k}(\xi)||\widehat{\nu_{k+s}}(\xi)|\le C \|\Omega\|_{L^q(\mathbf S^{n-1})} 2^{ -\gamma s}\min\{|2^k \xi|, |2^k \xi|^{-\gamma}\}$$ and
$$|1-\widehat{\phi_k}(\xi)||\widehat{\nu_{k+s}}(\xi)||\psi(2^{k-l}\xi)|\le C\|\Omega\|_{L^q(\mathbf S^{n-1})} 2^{ -\gamma s}  \min\{2^l, 2^{-\gamma l}\}.$$
Applying the above estimates and the Littlewood-Paley theory, we get
\begin{align}
\nonumber\|G_{s}^{l}f\|_{L^2}&\le C\|\Omega\|_{L^q(\mathbf S^{n-1})} 2^{ -\gamma s}  \min\{2^l, 2^{-\gamma l}\} \bigg\|\bigg(\dsum_{k\in \Bbb Z}|\Delta_{l-k}f|^2\bigg)^{1/2}\bigg\|_{L^2}\\
\label{Gsl22}&\le C\|\Omega\|_{L^q(\mathbf S^{n-1})} 2^{ -\gamma s}  \min\{2^l, 2^{-\gamma l}\}\|f\|_{L^2}.
\end{align}
Now we give the weighted $L^p$ estimate of $G_{ s}^lf$ for $l\in \Bbb Z$ and $s\in\Bbb N\cup \{0\}$.
Define
 $T_{s,k}f=[(\delta_0-\phi_k)\ast\nu_{k+s}]\ast f$. If we accept the following estimate for a moment
\begin{equation}\label{Tskfk}
\bigg\|
 \bigg(\dsum_{k\in \Bbb Z}
 | T_{s,k}f_k|^2\bigg)^{1/2}\bigg\|_{L^p(w)}\le C_{p,w}\|\Omega\|_{L^q(\mathbf S^{n-1})}\bigg\|
 \bigg(\dsum_{k\in \Bbb Z}
 | f_k|^2\bigg)^{1/2}\bigg\|_{L^p(w)},
\end{equation}
whenever $w$ and $p$ satisfy the conditions (i) or (ii) in Theorem \ref{thm:H}, then by \eqref{Tskfk} and the weighted Littlewood-Paley theory  (see \cite{K80}), we get
\begin{align}
\nonumber\|G_{ s}^lf\|_{L^p(w)}&=\bigg\|
 \bigg(\dsum_{k\in \Bbb Z}
 | T_{s,k}\Delta_{l-k}^2f|^2\bigg)^{1/2}\bigg\|_{L^p(w)}\\
 \nonumber&\le
C_{p,w}\|\Omega\|_{L^q(\mathbf S^{n-1})}\bigg\|
 \bigg(\dsum_{k\in \Bbb Z}
 |\Delta_{l-k}^2f|^2\bigg)^{1/2}\bigg\|_{L^p(w)}\\
 \label{Gslpp}&\le
 C_{p,w}\|\Omega\|_{L^q(\mathbf S^{n-1})}\|f\|_{L^p(w)}.
\end{align}
Repeating the same argument in \cite{D93} or \cite{Wa} and
using Stein and Weiss's interpolation theorem with change of measure between \eqref{Gsl22} and  \eqref{Gslpp},  we get for some $\theta_0,\beta_0\in(0,1)$
\begin{equation}\label{Gsllppd}
\|G_{ s}^lf\|_{L^p(w)}\le C_{p,w} 2^{-\beta_0 s}2^{-\theta_0|l|}
\|\Omega\|_{L^q(\mathbf S^{n-1})}\|f\|_{L^p(w)}.
\end{equation}
 Combining \eqref{lnt2}, \eqref{Gsl} and \eqref{Gsllppd}, we obtain
 \begin{align*}\|\sup_{\lambda>0}\lambda[N_\lambda(\mathscr T^{2}f)(x)]^{1/2}\|_{L^p(w)}&\le\dsum_{s\ge 0}\dsum_{l\in \Bbb Z}\|G_{ s}^lf\|_{L^p(w)}\\&\le  C_{p,w}\dsum_{s \ge 0}\dsum_{l\in \Bbb Z}2^{-\beta_0 s}2^{- \theta_0|l|}
\|\Omega\|_{L^q(\mathbf S^{n-1})}\|f\|_{L^p(w)}\\&\le
 C_{p,w}\|\Omega\|_{L^q(\mathbf S^{n-1})}\|f\|_{L^p(w)},
\end{align*}
whenever $w$ and $p$ satisfy the conditions (i) or (ii) in Theorem  \ref{thm:H}.

 Now we turn to the proof of \eqref{Tskfk}. First of all, we have the following estimate
 \begin{align}\label{vecTsk}
\bigg\|
 \bigg(\dsum_{k\in \Bbb Z}
 | T_{s,k}f_k|^2\bigg)^{1/2}\bigg\|_{L^p(w)}
 \le  C_{p,w}\bigg\|
 \bigg(\dsum_{k\in \Bbb Z}
 | \nu_{k+s}\ast f_k|^2\bigg)^{1/2}\bigg\|_{L^p(w)},
\end{align}
 which follows from
\begin{equation}\label{pkfk}
\bigg\|
\bigg(\dsum_{k\in \Bbb Z}
 | \phi_k\ast f_k|^2\bigg)^{1/2}\bigg\|_{L^p(w)}\le  C_{p,w}\bigg\|
 \bigg(\dsum_{k\in \Bbb Z}
 | f_k|^2\bigg)^{1/2}\bigg\|_{L^p(w)}
\end{equation} for $1< p<\infty$ and $w\in A_{p}$ (see \cite{AJ}). Secondly, we claim that for $w$ and $p$ satisfying the conditions (i) or (ii) in Theorem  \ref{thm:H}
\begin{equation}\label{vukk}
\bigg\|
 \bigg(\dsum_{k\in \Bbb Z}
 | \nu_{k+s}\ast f_k|^2\bigg)^{1/2}\bigg\|_{L^p(w)}\le  C_{p,w}\|\Omega\|_{L^q(\mathbf S^{n-1})}\bigg\|
 \bigg(\dsum_{k\in \Bbb Z}
 | f_k|^2\bigg)^{1/2}\bigg\|_{L^p(w)}.
\end{equation}
Then \eqref{Tskfk} is a consequence of \eqref{vecTsk} and \eqref{vukk}.

Now we turn to the proof of claim \eqref{vukk}. We only show that \eqref{vukk} holds when $w$ and $p$ satisfy the condition (i) in Theorem  \ref{thm:H}, the proof is similar for the condition (ii).  Obviously,
\begin{equation}\label{max cotr}
\big\|\sup_{k\in \Bbb Z}
  |\nu_{k+s}\ast f_k|\big\|_{L^p(w)}\le \big\|
M_{\Omega}\big(\sup_{k\in \Bbb Z}|f_k|\big)\big\|_{L^p(w)},
\end{equation}
where $M_\Omega$ denotes the rough maximal operator defined by
$$M_\Omega g(x)=\sup_{t>0}\frac1{t^n}\int_{|y|<t}|\Omega(y')g(x-y)|dy.$$
By \cite{D93} (see also \cite[p.\,106]{LDY}), if $w$ and $p$ satisfy the conditions (i) or (ii) in Theorem  \ref{thm:H}, then
 \begin{equation}\label{rough maxi wei bd}
 \|M_{\Omega} g\|_{L^p(w)}\le C_{p,w}\|\Omega\|_{L^q(\mathbf S^{n-1})}\|
 g\|_{L^p(w)}.
\end{equation}
Hence, from \eqref{max cotr} and \eqref{rough maxi wei bd}, if $w$ and $p$ satisfy the conditions (i) or (ii) in Theorem  \ref{thm:H}, then
\begin{equation}\label{bd of l-infty}
\big\|\sup_{k\in \Bbb Z}
  |\nu_{k+s}\ast f_k|\big\|_{L^p(w)}\le C_{p,w}\|\Omega\|_{L^q(\mathbf S^{n-1})}\big\|
\sup_{k\in \Bbb Z}|f_k|\big\|_{L^p(w)}.\end{equation}
Using \eqref{bd of l-infty} under the condition (ii) and the duality, we see that if  $w$ and $p$ satisfy the condition (i), then
\begin{equation}\label{bd of l-1}
\big\|\sum_{k\in \Bbb Z}
  |\nu_{k+s}\ast f_k|\big\|_{L^p(w)}\le C_{p,w}\|\Omega\|_{L^q(\mathbf S^{n-1})}\|
\sum_{k\in \Bbb Z}|f_k|\big\|_{L^p(w)}.\end{equation}
Interpolating between \eqref{bd of l-1} and \eqref{bd of l-infty} (under the condition (i)), we show that \eqref{vukk} holds if $w$ and $p$ satisfy the condition (i) in Theorem  \ref{thm:H}.
\medskip

\noindent{\bf Estimate of \eqref{lnti} for $i=3$.} Similarly,  we have the following pointwise estimate
\begin{align}
\nonumber\sup_{\lambda>0}\lambda[N_\lambda(\mathscr T^3f)(x)]^{1/2}
&\le\dsum_{s<0}\big(\dsum_{k\in \Bbb Z}\big| \phi_k\ast\nu_{k+s}\ast
f(x)\big|^2\big)^{1/2}\\
\label{Nldt3}&:=\dsum_{s<0}H_{s}f(x).
\end{align}
We are reduced to establish the estimate of $\|H_{s}f\|_{L^p(w)}$ as sharp as possible so that we are able to sum up over all
negative integers $s$. By  the Minkowski inequality, we get
\begin{align}
\nonumber\|H_{s}f\|_{L^p(w)}&=\bigg\|
 \bigg(\dsum_{k\in \Bbb Z}
 |\sum_{l\in \Bbb Z} \phi_k\ast\nu_{k+s}\ast\Delta_{l-k}^2 f|^2\bigg)^{1/2}\bigg\|_{L^p(w)}\\
 \nonumber&\le\dsum_{l\in \Bbb Z}\bigg\|
 \bigg(\dsum_{k\in \Bbb Z}
 |\phi_k\ast\nu_{k+s}\ast\Delta_{l-k}^2 f|^2\bigg)^{1/2}\bigg\|_{L^p(w)}\\
 \label{dhs}&:=\dsum_{l\in \Bbb Z}\|H_{s}^lf\|_{L^p(w)},
\end{align}
where the multipliers $\{\Delta_{l-k}\}$ were defined in Estimate of \eqref{lnti} for $i=2$.
We first prove a rapid decay estimate of $\|H_{s}^lf\|_{L^2(\mathbb R^n)}$ for $l\in \Bbb Z$ and $s<0$.
Since  $\hbox{supp}\ \widehat{\phi_k} \subset \{\xi :|2^k\xi|\le 4\}$, we have
$$|\widehat{\phi_k}(\xi)\widehat{\nu_{k+s}}(\xi)|\le C\|\Omega\|_{L^q(\mathbf S^{n-1})} 2^{ s} \min\{|2^k \xi|, |2^k \xi|^{-\gamma}\}$$
and
$$|\widehat{\phi_k}(\xi)\widehat{\nu_{k+s}}(\xi)\psi(2^{k-l}\xi)|\le C\|\Omega\|_{L^q(\mathbf S^{n-1})} 2^{s}  \min\{2^l, 2^{-\gamma l}\}.$$
Applying the above estimates and the Littlewood-Paley theory, we get
\begin{align}
\nonumber\|H_{s}^lf\|_{L^2}&\le C\|\Omega\|_{L^q(\mathbf S^{n-1})} 2^{s}  \min\{2^l, 2^{-\gamma l}\} \bigg\|\bigg(\dsum_{k\in \Bbb Z}|\Delta_{l-k}f|^2\bigg)^{1/2}\bigg\|_{L^2}\\
\label{Hsl22}&\le C\|\Omega\|_{L^q(\mathbf S^{n-1})} 2^{s}  \min\{2^l, 2^{-\gamma l}\} \|f\|_{L^2}.
\end{align}
Now we give the weighted $L^p$ norm of $H_{s}^lf$ for $l\in \Bbb Z$ and $s<0$.
By \eqref{pkfk}, \eqref{vukk} and  the weighted Littlewood-Paley theory  (see \cite{K80}), we get
\begin{align}
\nonumber\|H_{s}^lf\|_{L^p(w)}&\le
  C_{p,w}\|\Omega\|_{L^q(\mathbf S^{n-1})}\bigg\|
 \bigg(\dsum_{k\in \Bbb Z}
 |\Delta_{l-k}^2f|^2\bigg)^{1/2}\bigg\|_{L^p(w)}\\
 \label{Hslppw}&\le
 C_{p,w}\|\Omega\|_{L^q(\mathbf S^{n-1})}\|f\|_{L^p(w)},
 \end{align}
 whenever $w$ and $p$ satisfy the conditions (i) or (ii).
Repeating the same argument in \cite{D93} or \cite{Wa} and using Stein and Weiss's interpolation theorem with change of measure between \eqref{Hsl22} and \eqref{Hslppw},  we get for some $\beta_1, \theta_1\in(0,1)$,
\begin{align*}\|H_{s}^lf\|_{L^p(w)}&\le  C_{p,w} 2^{ \beta_1s}2^{- \theta_1|l|}
\|\Omega\|_{L^q(\mathbf S^{n-1})}\|f\|_{L^p(w)}.
\end{align*}
It follows from \eqref{Nldt3} and \eqref{dhs} that
$$\begin{array}{cl}\|\sup_{\lambda>0}\lambda[N_\lambda^d(\mathscr T^{3}f)(x)]^{1/2}\|_{L^p}&\le\dsum_{s< 0}\dsum_{l\in \Bbb Z}\|H_{ s}^lf\|_{L^p(w)}\\&\le  C_{p,w}\dsum_{s< 0}\dsum_{l\in \Bbb Z}2^{\beta_1 s}2^{- \theta_1|l|}
\|\Omega\|_{L^q(\mathbf S^{n-1})}\|f\|_{L^p(w)}\\&\le
 C_{p,w}\|\Omega\|_{L^q(\mathbf S^{n-1})}\|f\|_{L^p(w)},\end{array}$$whenever $w$ and $p$ satisfy the conditions (i) or (ii).
We therefore finish the estimate of \eqref{lnti} in the case $i=3$. \qed

\bigskip

At the end of this section, let us give the proof of Proposition \ref{pro:DNP}.

\noindent {\bf Proof of Proposition \ref{pro:DNP}}.
We first introduce some notations. For $j\in\mathbb Z$ and $\beta=(m_1,\cdots,m_n)\in\mathbb Z^n$,  we denote the dyadic cube $\prod_{k=1}^n(m_k2^j,(m_k+1)2^j]$ in $\mathbb R^n$ by $Q_\beta^j$, and the set of all dyadic cubes with sidelength $2^j$ by $\mathscr D_j$. The conditional expectation of a locally integrable $f$ with respect to $\mathscr D_j$ is given by
$$
\mathbb E_jf(x)=\sum_{Q\in \mathcal D_j}\frac1{|Q|}\int_{Q}f(y)dy\cdot\chi_{Q}(x)
$$
for all $j\in\mathbb Z$.
We also define the dyadic martingale difference operator $\mathbb D_j$ as
$\mathbb D_jf(x)=\mathbb E_{j}f(x)-\mathbb E_{j-1}f(x)$. Thus for $f\in L^p({\Bbb R^n})$, by the Lebesgue differential theorem we see that
\begin{equation}\label{Mar dec of f}
f(x)=-\sum_j\mathbb D_jf(x)\qquad \text{a. e.}\ x\in{\Bbb R^n}.
\end{equation}
We need to use a known extrapolation result:
\begin{lemma}{\rm (\cite[Corollary 1.2]{CMP})}\label{extrapolation}
Let $T$ be a sublinear operator such that $T: L^1(w)\rightarrow L^{1,\infty}(w)$ for all $w\in A_1$. Then
$\|Tf\|_{L^p(w)}\le C\|f\|_{L^p(w)}$ for $1<p<\infty$ and $w\in A_p$.
\end{lemma}
Note that $N_\lambda$ is subadditive, then
\begin{equation}\nonumber
N_\lambda(\mathscr Uf)\le N_{\lambda/2}(\mathscr Df)+N_{\lambda/2}(\mathscr Ef),
\end{equation}
where
$$
\mathscr Df=\{\phi_k\ast f-\mathbb E_kf\}_k \ \text{and}\ \ \mathscr Ef=\{\mathbb E_kf\}_k.
$$
It has been proved in \cite[p.8]{KZ14} that
$$
\|\sup_{\lambda>0}\lambda\sqrt{N_\lambda(\mathscr{E}f)}\|_{L^p(w)}\leq C_p\|f\|_{L^p(w)},\ 1<p<\infty,\ w\in A_p.
$$
On the other hand, we observe that
\begin{equation}\nonumber
\sup_{\lambda>0}\lambda\sqrt{N_\lambda(\mathscr{D}f)}\le\big(\sum_{k\in\mathbb Z}|\phi_k\ast f-\mathbb E_kf|^2\big)^{1/2}:=\mathfrak Sf.
\end{equation}
By Lemma \ref{extrapolation}, we just need to prove
\begin{equation}\label{ww11}
\sup_{\alpha>0}\alpha w(\{x:\mathfrak Sf(x)>\alpha\})\le C\|f\|_{L^1(w)},\ \ w\in A_1.
\end{equation}
\par
For any fixed $\alpha>0$, we perform the Calder\'{o}n-Zygmund decomposition of $f$ at height $\alpha$ using dyadic cubes, then there exists $\Lambda\subseteq\mathbb Z\times\mathbb Z^n$ such that the collection of dyadic cubes $\{Q_{\beta}^j\}_{(j,\beta)\in\Lambda}$ are disjoint and the following hold:
  \begin{itemize}
  \item[{\rm(i)}]
  $|\bigcup_{(j,\beta)\in\Lambda}Q_\beta^j|\le \frac1\alpha\|f\|_{L^1(\mathbb R^n)}$;
  \item[{\rm(ii)}]
  $|f(x)|\le \alpha$, if $x\not\in\bigcup_{(j,\beta)\in\Lambda}Q_\beta^j$;
   \item[{\rm(iii)}]
   $\alpha<\frac1{|Q_\beta^j|}\int_{Q_\beta^j}|f(x)|dx\le2^n\alpha$ for each $(j,\beta)\in\Lambda$.
   \end{itemize}
 We set
\begin{equation}\nonumber
 g(x)=\left\{
 \begin{array}{ll}
 f(x),& \text{if}\ x\not\in \bigcup_{(j,\beta)\in\Lambda}Q_{\beta}^j,\\
 \frac1{|Q_{\beta}^j|}\dint_{Q_{\beta}^j}f(y)dy,& \text{if}\ x\in Q_{\beta}^j,(j,\beta)\in\Lambda
 \end{array}
 \right.
 \end{equation}
 and
 \begin{equation}\nonumber
 b(x)=\sum_{(j,\beta)\in\Lambda}[f(x)-\mathbb E_jf(x)]\chi_{Q_{\beta}^j}(x):=\sum_{(j,\beta)\in\Lambda}b_{j,\beta}(x).
 \end{equation}
Clearly, $f=g+b$, $\|g\|_{L^\infty(\mathbb R^n)}\le 2\alpha$, $\|g\|_{L^1(\mathbb R^n)}\le\|f\|_{L^1(\mathbb R^n)}$ and $\|b\|_{L^1(\mathbb R^n)}\le 2\|f\|_{L^1(\mathbb R^n)}$.
\par
We accept the following fact for a moment
\begin{equation}\label{w22}
\|\mathfrak Sf\|_{L^2(w)}\le C\|f\|_{L^2(w)},\ \ w\in A_1,
\end{equation}
which will be proved later. By \eqref{w22}, the definition of $g$ and \eqref{eq def A1}, we have
\begin{align*}
w(\{x:\mathfrak Sg(x)>\alpha\})&\le \frac{C}{\alpha^2}\|\mathfrak Sg\|_{L^2(w)}^2\le\frac{C}{\alpha^2}\int_{\mathbb R^n}|g(x)|^2w(x)dx\le \frac{C}{\alpha}\int_{\mathbb R^n}|g(x)|w(x)dx\\
&\le \frac{C}{\alpha}\int_{(\bigcup Q_{\beta}^j)^c}|f(x)|w(x)dx+\frac{C}{\alpha}\sum_{(j,\beta)\in\Lambda}\int_{Q_{\beta}^j}
|f(y)|\frac{w(Q_{\beta}^j)}{|Q_{\beta}^j|}dy\\
&\le \frac{C}{\alpha}\int_{(\bigcup Q_{\beta}^j)^c}|f(x)|w(x)dx+\frac{C}{\alpha}\sum_{(j,\beta)\in\Lambda}\int_{Q_{\beta}^j}
|f(y)|w(y)dy\\
&\le \frac{C}{\alpha}\int_{\mathbb R^n}|f(x)|w(x)dx.
\end{align*}
Let $\tilde{Q}_\beta^j$ be the cube with center of $Q_\beta^j$ and $16$ times sidelength of $Q_\beta^j$. Observe that $w\in A_1$ has the doubling property,
\begin{align*}
w(\bigcup_{(j,\beta)\in\Lambda} \tilde{Q}_\beta^j)&\le C\sum_{(j,\beta)\in\Lambda}w(Q_\beta^j)\le C\sum_{(j,\beta)\in\Lambda}\frac{w(Q_\beta^j)}{|Q_\beta^j|}|Q_\beta^j|\\
&\le C\sum_{(j,\beta)\in\Lambda}\inf_{x\in Q_\beta^j}w(x)\frac1{\alpha}\int_{Q_\beta^j}|f(y)|dy\\
&\le C\sum_{(j,\beta)\in\Lambda}\frac1{\alpha}\int_{Q_\beta^j}|f(y)|w(y)dy\\
&\le \frac{C}{\alpha}\int_{\mathbb R^n}|f(y)|w(y)dy.
\end{align*}
Notice that $\mathbb E_kb_{j,\beta}(x)=0$ for $x\not\in Q_\beta^j$ (see \cite[p.6726]{JSW08}), then
\begin{align*}
\alpha w(\{(\bigcup \tilde{Q}_\beta^j)^c:\mathfrak Sb(x)>\alpha\})&\le\int_{(\bigcup \tilde{Q}_\beta^j)^c}\sum_{k\in\mathbb Z}|\phi_k\ast b(x)-\mathbb E_kb(x)|w(x)dx\\
&\le\sum_{(j,\beta)\in\Lambda}\sum_{k\in\mathbb Z}\int_{(\tilde{Q}_\beta^j)^c}|\phi_k\ast b_{j,\beta}(x)|w(x)dx.
\end{align*}
Recall that $\phi$ is a Schwartz function, so
$$
|\phi_k(x)|\le \frac{C2^{-kn}}{(1+|2^{-k}x|)^{n+1}}\ \ \text{and}\ \ |\partial_x\phi_k(x)|\le \frac{C2^{-k(n+1)}}{(1+|2^{-k}x|)^{n+2}}.
$$
For $k\le j$, by \eqref{eq def A1}, we obtain
\begin{align*}
\int_{(\tilde{Q}_\beta^j)^c}|\phi_k\ast b_{j,\beta}(x)|w(x)dx&\le \int_{Q_\beta^j}|b_{j,\beta}(y)|\int_{(\tilde{Q}_\beta^j)^c}|\phi_k(x-y)|w(x)dxdy\\
&\le \int_{Q_\beta^j}|b_{j,\beta}(y)|\int_{|x-y|\ge4\cdot2^j}\frac{2^{-nk}}{|2^{-k}(x-y)|^{n+1}}w(x)dxdy\\
&\le C2^{k-j} \int_{Q_\beta^j}|b_{j,\beta}(y)|M(w)(y)dy\\
&\le C2^{k-j} \big[\int_{Q_\beta^j}|f(y)|w(y)dy+\int_{Q_\beta^j}|f(z)|\frac{w(Q_\beta^j)}{|Q_\beta^j|}dz\big]\\
&\le C2^{k-j} \int_{Q_\beta^j}|f(y)|w(y)dy.
\end{align*}
For $k\ge j$, we use the fact $\int_{Q_\beta^j}b_{j,\beta}(y)dy=0$. Let $z^j_\beta$ be the center of $Q_\beta^j$, then
\begin{align*}
\int_{(\tilde{Q}_\beta^j)^c}|\phi_k\ast b_{j,\beta}(x)|w(x)dx&\le \int_{Q_\beta^j}|b_{j,\beta}(y)|\int_{(\tilde{Q}_\beta^j)^c}|\phi_k(x-y)-\phi_k(x-z_\beta^j)|w(x)dxdy\\
&\le 2^{j-k}\int_{Q_\beta^j}|b_{j,\beta}(y)|\int_{\mathbb R^n}\frac{2^{-kn}}{(1+2^{-k}|x-y|)^{n+2}}w(x)dxdy\\
&\le C2^{j-k} \int_{Q_\beta^j}|b_{j,\beta}(y)|M(w)(y)dy\\
&\le C2^{j-k} \int_{Q_\beta^j}|f(y)|w(y)dy.
\end{align*}
Combining the above estimates, we obtain
\begin{align*}
\alpha w(\{(\bigcup \tilde{Q}_\beta^j)^c:\mathfrak Sb(x)>\alpha\})
&\le C\sum_{(j,\beta)\in\Lambda}\sum_{k\in\mathbb Z}2^{-|j-k|} \int_{Q_\beta^j}|f(y)|w(y)dy\le C\|f\|_{L^1(w)}.
\end{align*}

Now we turn to the proof of \eqref{w22}, which can be showed in a similar way as Lemma 3.2 in \cite{JSW08}. Estimate \eqref{w22} is a consequence of the following fact: there exists a $\theta>0$ such that
\begin{equation}\label{de of di}
\|\phi_{k+j}\ast\mathbb D_jf-\mathbb E_{k+j}\mathbb D_{j}f\|_{L^2(w)}\le 2^{-\theta|k|}\|\mathbb D_jf\|_{L^2(w)}.
\end{equation}
 Indeed, using \eqref{Mar dec of f}, \eqref{de of di} and the Minkowski inequality,
\begin{align*}
\|\mathfrak Sf\|_{L^2(w)}&\le \big(\sum_k(\sum_j\|\phi_k\ast \mathbb D_jf(x)-\mathbb E_k\mathbb D_jf\|_{L^2(w)})^2\big)^{1/2}\\
&\le\big(\sum_k(\sum_j2^{-\theta|k-j|}\|\mathbb D_jf\|_{L^2(w)})^2\big)^{1/2}\\
&\le C_\theta(\sum_j\|\mathbb D_jf\|_{L^2(w)}^2)^{1/2}\le C_\theta\|f\|_{L^2(w)},
\end{align*}
the last inequality follows from the fact that the dyadic martingale square function is bounded on $L^2(w)$, see for instance \cite[Theorem 3.6]{Buc93} or \cite{LPR10}.
\par
Finally we prove \eqref{de of di}. When $k\ge0$, $\mathbb E_{k+j}\mathbb D_jf=0$. In \cite[p.6722]{JSW08}, Jones {\it et.al} have proved that $|\phi_{k+j}\ast\mathbb D_jf|\le C2^{-k}M(\mathbb D_jf)$. The weighted $L^p$-boundedness of the Hardy-Littlewood maximal function implies
\begin{equation}\label{de re kb0}
\|\phi_{k+j}\ast\mathbb D_jf-\mathbb E_{k+j}\mathbb D_{j}f\|_{L^2(w)}\le 2^{-k}\|\mathbb D_jf\|_{L^2(w)}.
\end{equation}
When $k<0$, $\mathbb E_{k+j}\mathbb D_jf=\mathbb D_jf$. Thus
\begin{align*}
\phi_{k+j}\ast\mathbb D_jf(x)-\mathbb E_{k+j}\mathbb D_{j}f(x)&=\int_{|y|\le 2^{k+j}}\phi_{k+j}(y)[\mathbb D_jf(x-y)-\mathbb D_jf(x)]dy\\
&\quad+\sum_{d\ge1}\int_{E_{k,j,d}}\phi_{k+j}(y)[\mathbb D_jf(x-y)-\mathbb D_jf(x)]dy\\
&:= I_0(x)+\sum_{d\ge1}I_d(x),
\end{align*}
where $E_{k,j,d}=\{y:2^{k+j+d-1}<|y|\le 2^{k+j+d+1}\}$. For a fixed integer $k\le -8$, we estimate $\|I_d\|_{L^2(w)}$ for the cases $d>|k|/2$ and $0\le d\le|k|/2$, respectively.
 \par
 For the case $d>|k|/2$, $y\in E_{k,j,d}$ and $N\ge n+1$, we have $|\phi_{k+j}(y)|\le C2^{-Nd}2^{-n(k+j)}$, since $\phi$ is a Schwartz function. Using the Minkowski inequality,
\begin{align*}
\|I_d\|_{L^2(w)}&\le C2^{-Nd-n(k+j)}\big[\int_{\mathbb R^n}\big(\int_{E_{k,j,d}}|\mathbb D_jf(x-y)-\mathbb D_jf(x)|dy\big)^2w(x)dx\big]^{1/2}\\
&\le C2^{-Nd-n(k+j)}\int_{E_{k,j,d}}\big[\int_{\mathbb R^n}(|\mathbb D_jf(x-y)|^2+|\mathbb D_jf(x)|^2)w(x)dx\big]^{1/2}dy\\
&\le C2^{-Nd-n(k+j)}\int_{E_{k,j,d}}\big[\int_{\mathbb R^n}|\mathbb D_jf(x-y)|^2w(x)dx\big]^{\frac12}dy+\frac{C}{2^d}\|\mathbb D_jf\|_{L^2(w)}.
\end{align*}
Let $\hat{Q}$ be the dyadic cube containing $Q$ with sidelength $2l(Q)$. We write
$$\mathbb D_jf=\sum_{Q\in\mathscr D_{j-1}}a_Q\chi_Q,\ \text{where}\ a_Q=\frac1{|\hat{Q}|}\int_{\hat{Q}}f(z)dz-\frac1{|Q|}\int_{Q}f(z)dz.$$
Note that $Q^k_\beta\bigcap Q^k_\alpha=\emptyset$ when $\beta\neq\alpha$. By a trivial calculation,
$$\|\mathbb D_jf\|_{L^2(w)}=\big(\sum_{Q\in\mathscr D_{j-1}}a_Q^2w(Q)\big)^{1/2}.$$
The Cauchy-Schwarz inequality and the properties of $A_1$ weight imply
\begin{align*}
&\int_{E_{k,j,d}}\big[\int_{\mathbb R^n}|\mathbb D_jf(x-y)|^2w(x)dx\big]^{\frac12}dy\\
&\le C\int_{E_{k,j,d}}\big(\sum_{Q\in\mathscr D_{j-1}}a_Q^2\int_{\mathbb R^n}\chi_Q(x-y)w(x)dx\big)^{1/2}dy\\
&\le C|E_{k,j,d}|^{1/2}\big(\sum_{Q\in\mathscr D_{j-1}}a_Q^2\int_{Q}\int_{E_{k,j,d}}w(x+y)dydx\big)^{1/2}\\
&\le C|E_{k,j,d}|^{1/2}2^{(k+j+d+1)n/2}\big(\sum_{Q\in\mathscr D_{j-1}}a_Q^2\int_{Q}M(w)(x)dx\big)^{1/2}\\
&\le C|E_{k,j,d}|(\sum_{Q\in\mathscr D_{j-1}}a_Q^2w(Q))^{1/2}\\
&\le C|E_{k,j,n}|\|\mathbb D_jf\|_{L^2(w)}.
\end{align*}
Then, we conclude that
\begin{align}\label{Inw2nl}
\|I_d\|_{L^2(w)}\le C\big(2^{-Nd-n(k+j)}|E_{k,j,d}|+2^{-d}\big)\|\mathbb D_jf\|_{L^2(w)}\le C2^{-d}\|\mathbb D_jf\|_{L^2(w)}.
\end{align}

For each $0\le d\le|k|/2$ and $Q_\beta^{j-1}$, we have $\mathbb D_jf(x-y)=\mathbb D_jf(x)$ for $y\in E_{k,j,d}$ and $x\in Q_\beta^{j-1}$ such that $\hbox{dist}(x,(Q_\beta^{j-1})^c)\ge 2^{j+d+k+1}$. By (1.14), there is a $\delta>0$ such that
\begin{align*}
\|I_d\|^2_{L^2(w)}&\le C\sum_{\beta}\sup_{x\in Q_\beta^{j-1}}|I_d(x)|^2w(\{x\in Q_\beta^{j-1}:dist(x,(Q_\beta^{j-1})^c)\le 2^{j+d+k+1}\})\\
&\le C\sum_{\beta}\frac{|\{x\in Q_\beta^{j-1}:dist(x,(Q_\beta^{j-1})^c)\le 2^{j+d+k+1}\}|^{\delta}}{|Q_\beta^{j-1}|^{\delta}}w(Q_\beta^{j-1})\sup_{x\in Q_\beta^{j-1}}|I_d(x)|^2\\
&\le C2^{-\delta|k|n/2}\sum_{\beta}w(Q_\beta^{j-1})\sup_{x\in Q_\beta^{j-1}}|I_d(x)|^2.
\end{align*}
For $x\in Q_\beta^{j-1}$, $|I_d(x)|\le C\sup_{y\in B(z_\beta^{j-1},2^{j+1})}|\mathbb D_jf(y)|$. Note that there are at most $9^n$ $Q^{j-1}_\alpha$'s jointing with $B(z_\beta^{j-1},2^{j+1})$. Therefore there exists a multi-index $\alpha_\beta$ such that
$$
\sup_{x\in B(z_\beta^{j-1},2^{j+1})}|\mathbb D_jf(x)|=|a_{Q^{j-1}_{\alpha_\beta}}|=|\mathbb D_jf(x)|\chi_{Q^{j-1}_{\alpha_\beta}}(x).$$
Thus we obtain
\begin{align}
\nonumber\|I_d\|^2_{L^2(w)}&\le C2^{-\delta|k|n/2}\sum_{\beta}w(Q_\beta^{j-1})|a_{Q^{j-1}_{\alpha_\beta}}|^2\\
\nonumber&\le C2^{-\delta|k|n/2}\sum_{\beta}\frac{w(Q_\beta^{j-1})}{w(B(z_\beta^{j-1},2^{j+1}))}\int_{B(z_\beta^{j-1},2^{j+1})}|\mathbb D_jf(x)|^2\chi_{Q^{j-1}_{\alpha_\beta}}(x)w(x)dx\\
\label{Inw2ns}&\le  C2^{-\delta|k|n/2}\|\mathbb D_jf\|_{L^2(w)}^2.
\end{align}
For $k\le -8$, by estimates \eqref{Inw2nl} and \eqref{Inw2ns}, we have
\begin{align}
\nonumber\|\phi_{k+j}\ast\mathbb D_jf-\mathbb E_{k+j}\mathbb D_{j}f\|_{L^2(w)}&\le \|\sum_{0\le d\le|k|/2}I_d\|_{L^2(w)}+\|\sum_{d>|k|/2}I_d\|_{L^2(w)}\\
\nonumber&\le C|k|2^{-\delta|k|n/4}\|\mathbb D_jf\|_{L^2(w)}+C\sum_{d\ge|k|/2}\frac1{2^d}\|\mathbb D_jf\|_{L^2(w)}\\
\label{de re kl0}&\le C2^{-\theta|k|}\|\mathbb D_jf\|_{L^2(w)}
\end{align}
for some $\theta>0$. Finally \eqref{de re kb0} and \eqref{de re kl0} together imply the desired estimate \eqref{de of di} and we finish the proof of Proposition \ref{pro:DNP}. \qed

\section{Proof of Theorem \ref{thm:H} (II)}

In this section we will finish the proof of \eqref{short variation singular}.
For $t\in[1,2]$, we define $\nu_{0,t}$ as
$$\nu_{0,t}(x)=\frac{\Omega(x')}{|x|^n}\chi_{_{\{t\leq|x|\leq2\}}}(x)$$
and $\nu_{j,t}(x)={2^{-jn}}\nu_{0,t}(2^{-j}x)$ for $j\in\mathbb{Z}$. Observe that
$V_{2,j}(\mathcal{T}f)(x)$ is just the strong $2$-variation function of the family
$\{\nu_{j,t}\ast f(x)\}_{t\in[1,2]}$, the Minkowski inequality implies
\begin{align}
S_{2}(\mathcal{T}f)(x)&=\Big(\sum_{j\in\mathbb{Z}}|V_{2,j}(\mathcal{T}f)(x)|^2\Big)^{\frac{1}{2}}
=\Big(\sum_{j\in\mathbb{Z}}\|\{\nu_{j,t}\ast f(x)\}_{t\in[1,2]}\|_{V_2}^2\Big)^{\frac{1}{2}}\nonumber\\
&\leq\sum_{k\in\mathbb{Z}}\Big(\sum_{j\in\mathbb{Z}}\|\{\nu_{j,t}\ast(\Delta^2_{k-j}
f)(x)\}_{t\in[1,2]}\|_{V_2}^2\Big)^{\frac{1}{2}}\nonumber\\
\label{Decom of S-2}&:=
\sum_{k\in\mathbb{Z}}S_{2,k}(\mathcal{T}f)(x).
\end{align}
By \eqref{Decom of S-2}, to get \eqref{short variation singular} it remains to show that
if $p$ and $w$ satisfy the conditions (i) or (ii) in Theorem \ref{thm:H}, then there exists an $\varepsilon>0$ such that for all $k\in\mathbb{Z}$
\begin{equation}\label{est of S-2k}
\|S_{2,k}(\mathcal{T}f)\|_{L^{p}(w)}\leq
 C_{p,w}2^{-\varepsilon|k|}\|\Omega\|_{L^q(\mathbf{S}^{n-1})}\|f\|_{L^{p}(w)}.
\end{equation}
Note that the authors of \cite{DHL} proved that for given $1<p<\infty,$ there exists a constant $\delta\in (0,1)$ such that for $k\in \Bbb Z$,
\begin{equation}\label{s2kp}
\|S_{2,k}(\mathcal{T}f)\|_{L^p}\leq
C_{p}2^{-\delta |k|}\|\Omega\|_{L^q(\mathbf{S}^{n-1})}\|f\|_{L^p}.
\end{equation}
Thus, if we can prove that when $p$ and $w$ satisfy the conditions (i) or (ii) in Theorem \ref{thm:H}, for all $k\in \Bbb Z$
\begin{equation}\label{s2kpw}
\|S_{2,k}(\mathcal{T}f)\|_{L^{p}(w)}\leq
 C_{p,w}\|\Omega\|_{L^q(\mathbf{S}^{n-1})}\|f\|_{L^{p}(w)},
\end{equation}
then by applying Stein and Weiss's interpolation theorem with change of measure between \eqref{s2kp}  and  \eqref{s2kpw}, there is a constant $0<\delta'<1$ such that
$$
\|S_{2,k}(\mathcal{T}f)\|_{L^{p}(w)}\leq
 C_{p,w}2^{-\delta'|k|}\|\Omega\|_{L^q(\mathbf{S}^{n-1})}\|f\|_{L^{p}(w)}.
$$
In other words, we get \eqref{est of S-2k} and complete the proof of \eqref{short variation singular}. Therefore, we reduce the proof of \eqref{short variation singular}  to showing \eqref{s2kpw}.
The estimate of \eqref{s2kpw} is dealt with differently according to whether $q\ge 2$ or $q<2$.

\subsection {The estimate of \eqref{s2kpw} for $q\ge 2$}\quad

We have the following observation:
$$\begin{array}{cl}S_{2,k}(\mathcal{T}f)(x)&=\Big(\dsum_{j\in\mathbb{Z}}\dsup_{\substack
{t_1<\cdots<t_N\\
[t_l,t_{l+1}]\subset[1,2]}}\dsum_{l=1}^{N-1}|\nu_{j,t_l}\ast\Delta_{k-j}^2
f(x)-\nu_{j,t_{l+1}}\ast\Delta_{k-j}^2
f(x)|^2\Big)^{\frac{1}{2}}\\&=
\bigg(\dsum_{j\in \mathbb{Z}}\sup_{\substack{t_{1}<\cdots<t_{N}\\ [t_{l},t_{l+1}]\subset[1,2]}}\dsum_{l=1}^{N-1}
\bigg| \int_{2^{j}t_{l}\leq|y|\leq2^{j+1}}\frac{\Omega(y')}{|y|^{n}}\Delta^{2}_{k-j}f(x-y)\,dy\\&\quad-\dint_{2^{j}t_{l+1}\leq| y|\leq2^{j+1}}\frac{\Omega(y')}{|y|^{n}}\Delta^{2}_{k-j}f(x-y)\,dy\bigg|^{2}\bigg)^{\frac{1}{2}}\\
&=\bigg(\dsum_{j\in \mathbb{Z}}\sup_{\substack{t_{1}<\cdots<t_{N}\\ [t_{l},t_{l+1}]\subset[1,2]}}\dsum_{l=1}^{N-1}
\bigg| \dint_{2^{j}t_{l}\leq|y|<2^{j}t_{l+1}}\frac{\Omega(y')}{|y|^{n}}\Delta^{2}_{k-j}f(x-y)\,dy\bigg|^{2}\bigg)^{\frac{1}{2}}\\
&\leq\bigg(\dsum_{j\in \mathbb{Z}}\sup_{\substack{t_{1}<\cdots<t_{N}\\ [t_{l},t_{l+1}]\subset[1,2]}}\dsum_{l=1}^{N-1}
\bigg( \dint_{2^{j}t_{l}\leq|y|<2^{j}t_{l+1}}\frac{|\Omega(y')|}{|y|^{n}}|\Delta^{2}_{k-j}f(x-y)|\,dy\bigg)^{2}\bigg)^{\frac{1}{2}}\\
&\leq\bigg(\dsum_{j\in \mathbb{Z}}\sup_{\substack{t_{1}<\cdots<t_{N}\\ [t_{l},t_{l+1}]\subset[1,2]}}\bigg(\dsum_{l=1}^{N-1}
\dint_{2^{j}t_{l}\leq|y|<2^{j}t_{l+1}}\frac{|\Omega(y')|}{|y|^{n}}|\Delta^{2}_{k-j}f(x-y)|\,dy\bigg)^{2}\bigg)^{\frac{1}{2}}.\end{array}$$ Let $B_l=\{y: 2^{j}t_{l}\leq|y|<2^{j}t_{l+1}\}.$ Since $t_1<\cdots<t_N,$ then $B_l's$ are disjoint and $\bigcup_{l=1}^{N-1}B_l=\{y:2^{j}t_{1}\leq|y|<2^{j}t_{N}\}$. We get
$$\begin{array}{cl}&\bigg(\dsum_{j\in \mathbb{Z}}\sup_{\substack{t_{1}<\cdots<t_{N}\\ [t_{l},t_{l+1}]\subset[1,2]}}\bigg(\dsum_{l=1}^{N-1}
\dint_{2^{j}t_{l}\leq|y|<2^{j}t_{l+1}}\frac{|\Omega(y')|}{|y|^{n}}|\Delta^{2}_{k-j}f(x-y)|\,dy\bigg)^{2}\bigg)^{\frac{1}{2}}\\
&=\bigg(\dsum_{j\in \mathbb{Z}}\sup_{\substack{t_{1}<t_{N}\\ [t_{1},t_{N}]\subset[1,2]}}\bigg(
\dint_{2^{j}t_{1}\leq|y|<2^{j}t_{N}}\frac{|\Omega(y')|}{|y|^{n}}|\Delta^{2}_{k-j}f(x-y)|\,dy\bigg)^{2}\bigg)^{\frac{1}{2}}\\
&=\bigg(\dsum_{j\in \mathbb{Z}}\bigg(
\dint_{2^{j}\leq|y|<2^{j+1}}\frac{|\Omega(y')|}{|y|^{n}}|\Delta^{2}_{k-j}f(x-y)|\,dy\bigg)^{2}\bigg)^{\frac{1}{2}}\\
&\le \bigg(\dsum_{j\in \mathbb{Z}}\bigg(\frac{1}{2^{jn}}
\dint_{|y|<2^{j+1}}|\Omega(y')||\Delta^{2}_{k-j}f(x-y)|\,dy\bigg)^{2}\bigg)^{\frac{1}{2}}\\&\le C\Big(\dsum_{j\in\mathbb{Z}}|M_{\Omega}(\Delta_{k-j}^2
f)(x)|^2\Big)^{\frac{1}{2}},\end{array}$$ where $M_\Omega$ is the rough maximal operator defined in Section 2.
Therefore, we get
 \begin{align}\label{s2km}
S_{2,k}(\mathcal{T}f)(x)\le C\Big(\dsum_{j\in\mathbb{Z}}|M_{\Omega}(\Delta_{k-j}^2
f)(x)|^2\Big)^{\frac{1}{2}}.
\end{align}
To continue the proof, we need the following proposition which will be proved later.
\begin{proposition}\label{m2pwp} For $q\ge2$, if $p$ and $w$ satisfy the conditions (i) or (ii) in Theorem \ref{thm:H}, then
 \begin{equation}\label{m2pw}
 \bigg\|\bigg(\sum_{j\in \Bbb Z}
|M_\Omega f_j|^2\bigg)^{1/2}\bigg\|_{L^p(w)}\le  C_{p,w}(1+\|\Omega\|_{L^q(\mathbf{S}^{n-1})})\bigg\|\bigg(\sum_{j\in \Bbb Z}
|f_j|^2\bigg)^{1/2}\bigg\|_{L^p(w)}.
\end{equation}
\end{proposition}
Using \eqref{s2km}, \eqref{m2pw} and the weighted Littlewood-Paley theory (see for instance \cite{K80}), we get for $q\ge2,$
$$\begin{array}{cl}\|S_{2,k}(\mathcal{T}f)\|_{L^p(w)}&\le   C_{p,w}\|\Omega\|_{L^q(\mathbf{S}^{n-1})}\bigg\|\Big(\dsum_{j\in\mathbb{Z}}|\Delta_{k-j}^2
f|^2\Big)^{\frac{1}{2}}\bigg\|_{L^p(w)}\le   C_{p,w}\|\Omega\|_{L^q(\mathbf{S}^{n-1})}\|f\|_{L^p(w)},\end{array}$$ whenever $w$ and $p$ satisfy the conditions (i) or (ii).

Now we prove  Proposition \ref{m2pwp}. Write
\begin{align*}
|\Omega(x')|&=\bigg[|\Omega(x')|-\frac1{\omega_{n-1}}\int_{\mathbf S^{n-1}}|\Omega(y')|d\sigma(y')\bigg]+\frac1{\omega_{n-1}}\int_{\mathbf S^{n-1}}|\Omega(y')|d\sigma(y')\\
&:=\Omega_0(x')+C(\Omega,n),
\end{align*}
where $\omega_{n-1}$ denotes the area of $\mathbf S^{n-1}$. It is easy to check that
\begin{equation}\label{rou max contrl}
{M}_\Omega f(x)\le
Cg_{\Omega_0}(|f|)(x)+CMf(x),
\end{equation}
where the operator $g_{\Omega_0}$ is defined by
$$
g_{\Omega_0}(f)(x)=\big(\sum_{k}|T_{k,\Omega_0}f(x)|^2\big)^{1/2}\quad \mbox{with}\quad T_{k,\Omega_0}f(x)=\int_{2^k<|y|\le2^{k+1}}\frac{\Omega_0(y')}{|y|^n}f(x-y)dy
$$
and $M$ denotes the Hardy-Littlewood maximal operator.
By the properties of $A_p$ weights (see \cite{GR}), if $w$ and $p$ satisfy the conditions (i) or (ii) in Theorem \ref{thm:H}, then $w\in A_p$. Thus, by \eqref{rou max contrl} and using the weighted norm inequality of the operator $M$ for vector-valued functions (see \cite[Theorem 3.1]{AJ}),
to get \eqref{m2pw} it suffices to prove that if $q\ge2$ and  $w$, $p$ satisfy the conditions (i) or (ii) in Theorem \ref{thm:H},
\begin{equation}\label{v-v ineq 0f G-fun}
 \Big\|\Big(\sum_{j\in \Bbb Z}
|g_{\Omega_0} (f_j)|^2\Big)^{1/2}\Big\|_{L^p(w)}\le C_{p,w}\|\Omega_0\|_{L^q(\mathbf S^{n-1})}\Big\|
\Big(\sum_{j\in \Bbb Z}|f_j|^2\Big)^{1/2}\Big\|_{L^p(w)}.
\end{equation}

 Consider the linear operators $T_{\epsilon, \Omega_0}f=\sum_k \epsilon_k T_{k,\Omega_0}f$, where $\epsilon=\{\epsilon_k\}$ is  a sequence with $\epsilon_k=\pm1$. By the usual argument with Rademacher function \cite{St}, to show \eqref{v-v ineq 0f G-fun} it needs only to  prove that if $q\ge2$ and $w$, $p$ satisfy the conditions (i) or (ii) in Theorem \ref{thm:H}, then
\begin{equation}\label{Te2pw}
\Big\|\Big(\sum_{j\in \Bbb Z}|T_{\epsilon,\Omega_0} f_j|^2\Big)^{1/2}\Big\|_{L^p(w)}\le  C_{p,w}\|\Omega_0\|_{L^q(\mathbf S^{n-1})}\Big\|
\Big(\sum_{j\in \Bbb Z}|f_j|^2\Big)^{1/2}\Big\|_{L^p(w)},
\end{equation}
 where $C_{p,w}$ is independent of $\epsilon$. We need the following lemma,  which is only an application of \cite[Lemma 9.5.4]{G}.

\begin{lemma}\label{gh}
$(a)$ Let $q'\le 2<p<\infty$ and $w\in A_{p/q'}.$ Then there exists a constant $C_1=C_1(n,p, [w]_{A_{p/q'}})$ such that for every nonnegative function $g$ in $L^{(p/2)'}(w)$, there is a function $G(g) $ such that
\begin{enumerate}[\upshape (i)]
\item
   $g\le G(g)$

\item $\|G(g)\|_{L^{(p/2)'}(w)}\le 2\|g\|_{L^{(p/2)'}(w)}$

\item $[G(g) w ]_{A_{2/q'}}\le C_1.$\end{enumerate}
$(b)$ Let $q'<p<2$ and $w\in A_{p/q'}.$ Then there exists a constant $C_2=C_2(n,p, [w]_{A_{p/q'}})$ such that for every nonnegative function $h$ in $L^{\frac{p}{2-p}}(w)$, there is a function $H(h) $ such that

 \begin{enumerate}[\upshape (i)] \item $h\le H(h)$

 \item $\|H(h)\|_{L^{\frac{p}{2-p}}(w)}\le 2\|h\|_{L^{\frac{p}{2-p}}(w)}$

  \item $[H(h)^{-1}
  w]_{A_{2/q'}} \le C_2.$\end{enumerate}
 Moreover, both constants $C_1(n,p,B)$ and $C_2(n,p,B)$ increase as $B$ increases.
\end{lemma}

By the duality, we only need to prove \eqref{Te2pw} when $w$ and $p$ satisfy the condition (i) in Theorem \ref{thm:H}. We first prove \eqref{Te2pw} for $q'< p<\infty$.

{\it Case 1. $q'<p=2.$}\quad Note that $T_{\epsilon,\Omega_0}$ is weighted $L^2$ bounded for $A_{2/q'}$ weight uniformly in $\epsilon$ (see for instance \cite{D93} or \cite{Wa}),
so \eqref{Te2pw} is an obvious consequence of the fact above since $w\in A_{2/q'}$ in this case.

 {\it Case 2. $q'\le  2< p<\infty.$}\quad   Note that
 \begin{align*}
 \big\|\big(\dsum_{j\in \Bbb Z}|T_{\epsilon,\Omega_0} f_j|^2\big)^{1/2}\big\|_{L^p(w)}
=\sup\limits_{\|g\|_{ L^{(p/2)'}(w)}\le 1} \bigg|\dint_{{\Bbb R}^n}\dsum_{j\in \Bbb Z}
|T_{\epsilon,\Omega_0} f_j(x)|^2
g(x)w(x)\,dx\bigg|^{1/2}.
\end{align*}
 Since $w\in A_{p/q'}$ and $g\in L^{(p/2)'}(w)$, by Lemma \ref{gh} (a), there exists a function $G(|g|)$ such that
 \begin{align*}
 \big\|\big(\dsum_{j\in \Bbb Z}|T_{\epsilon,\Omega_0} f_j|^2\big)^{1/2}\big\|_{L^p(w)}\le \sup\limits_{\|g\|_{ L^{(p/2)'}(w)}\le 1} \bigg|\dsum_{j\in \Bbb Z}\dint_{{\Bbb R}^n}
|T_{\epsilon,\Omega_0} f_j(x)|^2\
G(|g|)(x)w(x)\,dx\bigg|^{1/2}.
\end{align*}
Note that $G(|g|)w\in A_{2/q'}$, hence by the weighted uniform $L^2$ boundedness of  $T_{\epsilon,\Omega_0}$, H\"{o}lder's inequality and Lemma \ref{gh} (a), we get
 \begin{align}
 \big\|\big(\dsum_{j\in \Bbb Z}|T_{\epsilon,\Omega_0} f_j|^2\big)^{1/2}\big\|_{L^p(w)}\notag
 &\le C\|\Omega_0\|_{L^q(\mathbf{S}^{n-1})}\sup\limits_{\|g\|_{ L^{(p/2)'}(w)}\le 1}\bigg|\dint_{{\Bbb R}^n}\dsum_{j\in \Bbb Z}
|f_j(x)|^{2}
G(|g|)(x)w(x)\,dx\bigg|^{1/2}\notag\\
&\le C\|\Omega_0\|_{L^q(\mathbf{S}^{n-1})} \sup\limits_{\|g\|_{ L^{(p/2)'}(w)}\le 1} \big\|\dsum_{j\in \Bbb Z}
|f_j|^{2}\big\|_{L^{p/2}(w)}^{1/2}\|G(|g|)\|_{L^{(\frac{p}{2})^{\prime}}(w)}^{1/2}\notag\\
&\le C\|\Omega_0\|_{L^q(\mathbf{S}^{n-1})}\sup\limits_{\|g\|_{ L^{(p/2)'}(w)}\le 1} \big\|\big(\dsum_{j\in \Bbb Z}
|f_j|^{2}\big)^{1/2}\big\|_{L^{p}(w)}\|g\|_{L^{(\frac{p}{2})^{\prime}}(w)}^{1/2}\notag\\
\label{Te2pw1}&\le C\|\Omega_0\|_{L^q(\mathbf{S}^{n-1})}\big\|\big(\dsum_{j\in \Bbb Z}
|f_j|^{2}\big)^{1/2}\big\|_{L^{p}(w)}.
\end{align}

{\it Case 3. $q'<p<2.$}\quad Since $w\in A_{p/q'}$ and $\big(\sum_{j\in \Bbb Z}|f_j|^2\big
)^{1/2}\in L^p(w)$, thus $h=\big(\sum_{j\in \Bbb Z}|f_j|^2\big)^{\frac{2-p}{2}}\in L^{\frac{p}{2-p}}(w).$
 Let $H(h)$ be the function in Lemma \ref{gh} (b). Then applying  H\"{o}lder's inequality, we have
\begin{align*}\big\|\big(\dsum_{j\in \Bbb Z}|T_{\epsilon,\Omega_0} f_j|^2\big)^{1/2}\big\|_{L^p(w)}
&=\big\|\dsum_{j\in \Bbb Z}
|T_{\epsilon,\Omega_0} f_j|^2H(h)^{-1} H(h)\big\|_{L^{\frac{p}{2}}(w)}^{1/2}\notag\\
&\le  \big\|\sum_{j\in \Bbb Z}
|T_{\epsilon,\Omega_0} f_j|^2H(h)^{-1}\big\|_{L^1(w)}^{1/2}\big\| H(h)\big\|_{L^{\frac{p}{2-p}}(w)}^{1/2}.
\end{align*}
  Since $H(h)^{-1}w\in A_{2/q'}$, using again the weighted $L^2$ uniform estimate of $T_{\epsilon,\Omega_0}$ and Lemma \ref{gh} (b), we have
\begin{align}
&\big\|\big(\dsum_{j\in \Bbb Z}|T_{\epsilon,\Omega_0} f_j|^2\big)^{1/2}\big\|_{L^p(w)}\notag\\
&\le  C\|\Omega_0\|_{L^q(\mathbf{S}^{n-1})}\bigg(\dsum_{j\in \Bbb Z}\dint_{{\Bbb R}^n}
|f_j(x)|^{2}H(h)(x)^{-1}w(x)\,dx\bigg)
^{1/2}\bigg\|\bigg(\dsum_{i\in \Bbb Z
} |f_i|^2\bigg)^{1-\frac{p}{2}}
\bigg\|_{L^{\frac{p}{2-p}}(w)}^{1/2}\notag\\
&\le C\|\Omega_0\|_{L^q(\mathbf{S}^{n-1})} \bigg(\dint_{{\Bbb R}^n}
\dsum_{j\in \Bbb Z}|f_j|^{2}
\big(\dsum_{i\in \Bbb Z} |f_i|^2\big)^{\frac{p}{2}-1}
w(x)\,dx\bigg)
^{1/2}\bigg\|\bigg(\dsum_{i\in \Bbb Z
} |f_i|^2\bigg)^{1/2}
\bigg\|_{L^{p}(w)}^{1-\frac{p}{2}}\notag\\
\label{Te2pw2}&=C\|\Omega_0\|_{L^q(\mathbf{S}^{n-1})}\bigg\|\bigg(\dsum_{j\in \Bbb Z}
|f_j|^{2}
\bigg)^{1/2}\bigg\|_{L^{p}(w)}.
\end{align}
Combining \eqref{Te2pw1} and \eqref{Te2pw2}, we show \eqref{Te2pw} for $w\in A_{p/q'}$ and $q'<p<\infty$.

Finally, the estimate of \eqref{Te2pw} for the endpoint case $p=q'$, $p\neq 1$ can follow by interpolating between \eqref{Te2pw} for some $p_0>q'$ and the unweighted estimate for some $p_1<q'$ (see for instance \cite{D93}) by using the properties of $A_p$ weight.

\begin{remark}\label{one key remark}
The main idea behind the proof of estimate (\ref{Te2pw}) is Rubio de Francia's extrapolation theorem \cite[Theorem 3]{Ru}. In particular, the result in the case $q>2$ and $p>q'$ is just a consequence of Rubio de Francia's extrapolation theorem by the fact that  $T_{\epsilon,\Omega_0}$ is bounded on $L^2(w)$ with $w\in A_{2/q'}$ (see \cite{D93} or \cite{Wa}). However, other cases $q=2$ or $p=q'$ can not be directly deduced from Theorem 3 in \cite{Ru}. Thus we prefer having given a detailed proof of estimate (\ref{Te2pw}) above.
\end{remark}

\subsection {The estimate of \eqref{s2kpw} for $q<2$}\quad

 For $q<2,$ we first give the  estimate of \eqref{s2kpw} for $p$ and $w$ satisfy the condition (ii) in Theorem \ref{thm:H}, that is,  $q<2,$ $1<p\le q$  and $ w^{-\frac{1}{(p-1)}}\in A_{p'/q'}$. We claim  that \eqref{m2pw} holds also for $q<2$, $p$ and $w$ satisfy the condition (ii) in Theorem \ref{thm:H}. In fact,
 by the weighted $L^p$ boundedness of $M_\Omega$ (see \cite{D93} or \cite{LDY}), we have
 $$
\bigg\|
 \bigg(\dsum_{j\in \Bbb Z}
 |M_\Omega f_j|^p\bigg)^{1/p}\bigg\|_{L^p(w)}\le C_{p,w}\|\Omega\|_{L^q(\mathbf S^{n-1})}\bigg\|
 \bigg(\dsum_{j\in \Bbb Z}
 | f_j|^p\bigg)^{1/p}\bigg\|_{L^p(w)}
$$
and
$$\big\|\sup_{j\in \Bbb Z}
  |M_\Omega f_j|\big\|_{L^p(w)}\le C_{p,w}\|\Omega\|_{L^q(\mathbf S^{n-1})}\|
\sup_{j\in \Bbb Z}|f_j|\big\|_{L^p(w)}
$$
for $1< p\le q$  and $ w^{-\frac{1}{(p-1)}}\in A_{p'/q'}$.
Interpolating between the above two estimates, we get
$$
\bigg\|
 \bigg(\dsum_{j\in \Bbb Z}
 |M_\Omega f_j|^2\bigg)^{1/2}\bigg\|_{L^p(w)}\le C_{p,w}\|\Omega\|_{L^q(\mathbf S^{n-1})}\bigg\|
 \bigg(\dsum_{j\in \Bbb Z}
 | f_j|^2\bigg)^{1/2}\bigg\|_{L^p(w)}
$$
for $q<2$, $1<p\le q$  and $ w^{-\frac{1}{(p-1)}}\in A_{p'/q'}$.
Thus, by  \eqref{s2km}, \eqref{m2pw}  and  the weighted Littlewood-Paley theory (see \cite{K80}), we
have
\begin{align}
\nonumber\|S_{2,k}(\mathcal{T}f)\|_{L^p(w)}&\le\bigg\|\Big(\dsum_{j\in\mathbb{Z}}|M_{\Omega}(\Delta_{k-j}^2
f)|^2\Big)^{\frac{1}{2}}\bigg\|_{L^p(w)}\\
\nonumber&\le C_{p,w}\|\Omega\|_{L^q(\mathbf{S}^{n-1})}\bigg\|\Big(\dsum_{j\in\mathbb{Z}}|\Delta_{k-j}^2
f|^2\Big)^{\frac{1}{2}}\bigg\|_{L^p(w)}\\
\label{s2kpwii}&\le  C_{p,w}\|\Omega\|_{L^q(\mathbf{S}^{n-1})}\|f\|_{L^p(w)},
\end{align}whenever  $q<2$,  $1<p\le q$  and $ w^{-\frac{1}{(p-1)}}\in A_{p'/q'}$.

Secondly, for $q<2,$ we estimate \eqref{s2kpw} for $p$ and $w$ satisfying the condition (i) in Theorem \ref{thm:H}, that is, $q<2,$ $q'\le p<\infty$ and $w\in A_{p/q'}$.
Using the idea of proving Lemma B.1 in \cite{KZ14}, we may show that for $x\in\Bbb R^n$,  \begin{equation}\label{cotr inq}
\|\mathfrak{a}_t\|_{\tilde{V}_2}\le8\|\mathfrak{a}_t\|_{X}^{1/2}
\Big\|\frac d{dt}\mathfrak{a}_t\Big\|_{X}^{1/2},
\end{equation}
where $X=L^2([1,2],\frac {dt}t)$ and $\mathfrak{a}_t=T_{j,t}
\Delta_{k-j}^2f(x)$ with $T_{j,t}h(x)=\nu_{j,t}\ast h(x)$. \eqref{cotr inq} means that
\begin{align*}
[S_{2,k}(\mathcal{T}f)(x)]^2\leq8\sum_{j\in\mathbb{Z}}\bigg(\dint_{1}^{2}|T_{j,t}
\Delta_{k-j}^2f(x)|^2\frac{dt}{t}\bigg)^{1/2}\bigg(\dint_{1}^{2}| \frac{d}{dt}T_{j,t}
\Delta_{k-j}^2f(x)|^2\frac{dt}{t}\bigg)^{1/2}.
\end{align*}
By the Cauchy-Schwarz inequality, we have for  any $1<p<\infty,$
\begin{align*}\nonumber
&\big\|S_{2,k}(\mathcal{T}f)\big\|^2_{L^p(w)}\\&\le
\bigg\|\bigg(\sum_{j\in\mathbb{Z}}\dint_{1}^{2}|T_{j,t}
\Delta_{k-j}^2f|^2\frac{dt}{t}\bigg)^{1/2}\bigg\|_{L^p(w)}\bigg\|\bigg(\sum_{j\in\mathbb{Z}}\dint_{1}^{2}|\frac{d}{dt}T_{j,t}
\Delta_{k-j}^2f|^2\frac{dt}{t}\bigg)^{1/2}\bigg\|_{L^p(w)}\\&:=\|I_{1,k}f\|_{L^p(w)}\cdot\|I_{2,k}f\|_{L^p(w)}.
\end{align*}
We now estimate $\|I_{1,k}f\|_{L^p(w)}$ and $\|I_{2,k}f\|_{L^p(w)},$ respectively.
For $\|I_{1,k}f\|_{L^p(w)}$, since  $p\ge q'>2$, by the Minkowski inequality, we get $$\begin{array}{cl}\|I_{1,k}f\|_{L^p(w)}&\le \bigg(\dint_{1}^{2}\bigg\|\bigg(\sum_{j\in\mathbb{Z}}|T_{j,t}
\Delta_{k-j}^2f|^2\bigg)^{1/2}\Big\|_{L^p(w)}^2\frac{dt}{t}\bigg)^{1/2}.\end{array}$$
By \cite{D93}, for $q'\le p<\infty$, $p\neq 1$ and $w\in A_{p/q'}$
   $$\big\|\sup_{j\in \Bbb Z,\,t\in [1,2)}|T_{j,t}f|\big\|_{L^p(w)}\le C\|M_\Omega f\|_{L^p(w)}\le C_{p,w}\|\Omega\|_{L^q(\mathbf{S}^{n-1})}\|f\|_{L^p(w)}.$$
Then similarly to the proof of   \eqref{Tskfk} and applying the weighted Littlewood-Paley theory with $w\in A_p$ (\cite{K80}), we get
$$\begin{array}{cl}\|I_{1,k}f\|_{L^p(w)}&\le C_{p,w}\|\Omega\|_{L^q(\mathbf{S}^{n-1})}\bigg(\dint_{1}^{2}\bigg\|\bigg(\sum_{j\in\mathbb{Z}}|
\Delta_{k-j}^2f|^2\bigg)^{1/2}\Big\|_{L^p(w)}^2\frac{dt}{t}\bigg)^{1/2}\le C_{p,w}\|\Omega\|_{L^q(\mathbf{S}^{n-1})}\|f\|_{L^p(w)}.\end{array}$$

Now we consider $\|I_{2,k}f\|_{L^p(w)}$. For any given Schwartz function $h$,
by the spherical coordinate transformation, a trivial calculation shows
\begin{align*}
\frac{d}{dt}T_{j,t}h(x)=\frac{d}{dt}[\nu_{j,t}\ast h(x)]&=\frac{d}{dt}\bigg[\int_{2^jt<|y|\le 2^{j+1}}\frac{\Omega(y')}{|y|^n}h(x-y)dy\bigg]\\&=\frac{d}{dt}\bigg[\int_{\mathbf S^{n-1}}\Omega(y')\int_{2^jt}^{2^{j+1}}\frac{1}{r}h(x-ry')drd\sigma(y')\bigg]\\
&=-\frac{1}{t}\int_{\mathbf S^{n-1}}\Omega(y')h(x-2^jty')d\sigma(y').
\end{align*}
For $t\in [1,2)$  and $\{h_j\}\in  L^2(\ell^2)({\Bbb R}^n), $ we define $$T_{j,t}^*h_j(x)=\int_{\mathbf S^{n-1}}|\Omega(y')||h_j(x-2^jty')|d\sigma(y').$$
It is easy to verify that$$\Big|\frac{d}{dt}T_{j,t}h_j(x)\Big|\le CT_{j,t}^*h_j(x).$$
 Next, we claim that for $q<2$, $q'\le p<\infty$ and $w\in A_{p/q'}$,
\begin{equation}\label{Tjt2pw}
\bigg\|\bigg(\dint_{1}^{2}\dsum_{j\in
\Bbb Z}\big|T_{j,t}^*
h_j\big|^2\frac{dt}{t}\bigg)^{1/2}\bigg\|_{L^p(w)}\le C_{p,w}\|\Omega\|_{L^q(\mathbf S^{n-1})}\bigg\|\bigg(\dsum_{j\in \Bbb Z}|h_j|^2\bigg)^{1/2}\bigg\|_{L^p(w)}.
\end{equation}
In fact,
for $t\in [1,2),$\begin{align*}
\big|T_{j,t}^*h_j(x)\big|^2\le \|\Omega\|_{L^q(\mathbf S^{n-1})}\int_{\mathbf S^{n-1}}|\Omega^{2-q}(y')|h_j(x-2^jty')|^2\,d\sigma(y').
\end{align*}
 Then for $q<2$ and $q'\le p<\infty,$ there exists a function $g\in L^{(\frac{p}{2})'}(w)$ with $\|g\|_{L^{(\frac{p}{2})'}(w)}=1,$ such that\begin{align*}
 &\bigg\|\bigg(\dint_{1}^{2}\dsum_{j\in
\Bbb Z}\big|T_{j,t}^*
h_j\big|^2\frac{dt}{t}\bigg)^{1/2}\bigg\|_{L^p(w)}^2\\&=\bigg|\int_{{\Bbb R}^n} \dint_{1}^{2}\dsum_{j\in
\Bbb Z}\big|T_{j,t}^*
h_j(x)\big|^2\frac{dt}{t}g(x)w(x)\,dx\bigg|\\&\le C\|\Omega\|_{L^q(\mathbf S^{n-1})}^q\int_{{\Bbb R}^n} \dsum_{j\in
\Bbb Z}\dint_{1}^{2}\int_{\mathbf S^{n-1}}|\Omega^{2-q}(y')||h_j(x-2^jty')|^2\,d\sigma(y')\frac{dt}{t}|g(x)|w(x)\,dx\\&\le C \|\Omega\|_{L^q(\mathbf S^{n-1})}^q\int_{{\Bbb R}^n} \dsum_{j\in
\Bbb Z}\int_{2^j<|y|\le 2^{j+1}}\frac{|\Omega^{2-q}(y')|}{|y|^n}|h_j(x-y)|^2\,dy|g(x)|w(x)\,dx\\&\le C \|\Omega\|_{L^q(\mathbf S^{n-1})}^q\int_{{\Bbb R}^n} M_{\Omega^{2-q}}(gw)(x)\dsum_{j\in
\Bbb Z}|h_j(x)|^2\,dx\\&\le C \|\Omega\|_{L^q(\mathbf S^{n-1})}^q\|M_{\Omega^{2-q}}(gw)\|_{L^{(\frac{p}{2})'}(w^{1/(1-p/2)})}\big\|\dsum_{j\in
\Bbb Z}|h_j|^2\big\|_{L^{\frac{p}{2}}(w)}\\&\le C_{p,w} \|\Omega\|_{L^q(\mathbf S^{n-1})}^{2}\bigg\|\bigg(\dsum_{j\in
\Bbb Z}|h_j|^2\bigg)^{1/2}\bigg\|_{L^{p}(w)}^2,
\end{align*}where in the above inequality we have used $$\bigg(\dint |M_{\Omega^{2-q}}(uw)(x)|^{(p/2)'}w(x)^{1/(1-p/2)}\,dx\bigg)^{1/(p/2)'}\le C_{p,w}\|\Omega\|_{L^q(\mathbf S^{n-1})}^{2-q}$$ for $\Omega^{2-q}\in L^{\frac{q}{2-q}}$ and $w\in A_{p/q'}$ (see \cite{D93} or \cite{LDY}).

Therefore, by (\ref{Tjt2pw}) and the weighted  Littlewood-Paley theory with $w\in A_p$ ( see \cite{K80}), we
have that for $q<2$ and $q'\le p<\infty$ and  $w\in A_{p/q'},$
$$\begin{array}{cl}\|I_{2,k}f\|_{L^p(w)}&\le C\bigg\|\bigg(\dint_{1}^{2} \dsum_{j\in
\Bbb Z}\big|T_{j,t}^*
\Delta_{k-j}f\big|^2
\frac{dt}{t}\bigg)^{1/2}\bigg\|_{L^p(w)}\\&\le C_{p,w}\|\Omega\|_{L^q(\mathbf S^{n-1})}\bigg\| \bigg(\dsum_{j\in
\Bbb Z}\big|
\Delta_{k-j}f\big|^2
\bigg)^{1/2}\bigg\|_{L^p(w)}\\&\le C_{p,w}\|\Omega\|_{L^q(\mathbf S^{n-1})}\|f\|_{L^p(w)}.\end{array}$$
Combined with the estimate of $\|I_{1,k}f\|_{L^p(w)}$,  we get  for $k\in \Bbb Z,$
\begin{align*}\nonumber
&\big\|S_{2,k}(\mathcal{T}f)\big\|_{L^p(w)}\le
C_{p,w}\|\Omega\|_{L^q(\mathbf S^{n-1})}\|f\|_{ L^p(w)},
\end{align*} whenever  $q<2$, $q'\le p<\infty,$ and $w\in A_{p/q'}.$\qed

\section{Proof of Theorem \ref{r o MO}}

We write
\begin{align*}
\Omega(x')&=\bigg[\Omega(x')-\frac1{\omega_{n-1}}\int_{\mathbf S^{n-1}}\Omega(y')d\sigma(y')\bigg]+\frac1{\omega_{n-1}}\int_{\mathbf S^{n-1}}\Omega(y')d\sigma(y'):=\Omega_0(x')+C(\Omega,n),
\end{align*}
where $\omega_{n-1}$ denotes the area of $\mathbf S^{n-1}$. Thus,
\begin{align*}
M_{\Omega,t}f(x)&=\frac1{t^n}\int_{|y|<t}\Omega_0(y')f(x-y)dy+C(\Omega,n)\frac1{t^n}\int_{|y|<t}f(x-y)dy\\&:=M_{\Omega_0,t}f+C(\Omega,n)M_tf,
\end{align*}
where $\Omega_0$ satisfies the cancelation condition \eqref{can of O}. Denote the operator family $\{M_{\Omega_0,t}\}_{t>0}$ by $\mathcal M_{\Omega_0}$ and $\{M_{t}\}_{t>0}$ by $\mathcal M$. It has been proved in \cite{KZ14} that $$\|\sup_{\lambda>0}\lambda\sqrt{N_{\lambda}(\mathcal Mf)}\|_{L^p(w)}\le C_{p,w}\|f\|_{L^p(w)}.$$
 To prove Theorem \ref{r o MO}, it suffices to show
\begin{equation}\label{nmopw}
\|\sup_{\lambda>0}\lambda\sqrt{N_{\lambda}(\mathcal M_{\Omega_0}f)}\|_{L^p(w)}\le C_{p,w}\|\Omega_0\|_{L^q(\mathbf S^{n-1})}\|f\|_{L^p(w)}.
\end{equation}
Thus by \eqref{lem:convert lemma}, the proof of  \eqref{nmopw} is reduced to prove
\begin{equation}\label{ndmopw}
\|\sup_{\lambda>0}\lambda\sqrt{N_{\lambda}(\{M_{\Omega_0, 2^k}f\}_{k\in \Bbb Z})}\|_{L^p(w)}\le
C_{p,w}\|\Omega_0\|_{L^q(\mathbf S^{n-1})}\|f\|_{L^p(w)}
\end{equation}
and
\begin{equation}\label{s2mopw}
\|S_2(\mathcal M_{\Omega_0}f)\|_{L^p(w)}\le
C_{p,w}\|\Omega_0\|_{L^q(\mathbf S^{n-1})}\|f\|_{L^p(w)}.
\end{equation}
For \eqref{ndmopw}, we define  $\sigma_{k}(y)={2^{-kn}}\Omega_0(y')\chi_{\{|y|<2^k\}}(y)$ for $k\in \Bbb Z$. For $q$, $w$ and $p$ satisfying the conditions (i) or (ii), we have
\begin{align*}
\|\sup_{\lambda>0}\lambda\sqrt{N_{\lambda}(\{M_{\Omega_0, 2^k}f\}_{k\in \Bbb Z})}\|_{L^p(w)}&\leq C\|\big(\sum_{k\in \Bbb Z}|f\ast\sigma_{k}|^2\big)^{1/2}\|_{L^p(w)}\leq C_{p,w}\|\Omega_0\|_{L^q(\mathbf{S}^{n-1})}\|f\|_{L^p(w)},
\end{align*}
where the second inequality is a known result in \cite{D93}.

For \eqref{s2mopw},  we define $\mu_{j,t}(x)={(2^{j}t)}^{-n}\Omega_0(x')\chi_{\{|x|<2^jt\}}(x)$ for $j\in\mathbb{Z}$ and $t\in[1,2)$. Observe that
$V_{2,j}(\mathcal{M}_{\Omega_0}f)(x)$ is just the strong  $2$-variation function of the family
$\{\mu_{j,t}\ast f(x)\}_{t\in[1,2]}$, hence
\begin{align*}
S_{2}(\mathcal{M}_{\Omega_0}f)(x)&=\Big(\sum_{j\in\mathbb{Z}}|V_{2,j}(\mathcal{M}_{\Omega_0}f)(x)|^2\Big)^{\frac{1}{2}}.
\end{align*}
 Similar to the proof of \eqref{short variation singular} in Theorem \ref{thm:H}, we get
 \begin{align*}
\|S_{2}(\mathcal{M}_{\Omega_0}f)\|_{L^p(w)}\le C_{p,w}\|\Omega_0\|_{L^q(\mathbf{S}^{n-1})}\|f\|_{L^p(w)},
\end{align*}
whenever $q$, $w$ and $p$ satisfy the conditions (i) or (ii). \qed

\vskip 1cm

 {\bf Acknowledgement.}\quad The authors would like to
express their deep gratitude to the referee for giving many valuable suggestions.

\bibliographystyle{amsplain}

\end{document}